\documentstyle[amsmath,amssymb
]{article}

\newcommand{\rom}[1]{{\rm #1}}
\newcommand{\wick}[1]{{:}\omega^{\otimes #1}{:}_\lambda}

\allowdisplaybreaks

\makeatletter\@addtoreset{equation}{section}\makeatother

\begin{document}

\setcounter{page}{1} \setcounter{section}{0} \thispagestyle{empty}

\newtheorem{definition}{Definition}[section]
\newtheorem{remark}{Remark}[section]
\newtheorem{proposition}{Proposition}[section]
\newtheorem{theorem}{Theorem}[section]
\newtheorem{corollary}{Corollary}[section]
\newtheorem{lemma}{Lemma}[section]

\newcommand{\indlim}{\operatornamewithlimits{ind\,lim}}
\newcommand{\Ffin}{{\cal F}_{\mathrm fin}}
\newcommand{\Fext}{{\cal F}_{\mathrm ext}}
\newcommand{\D}{{\cal D}}
\newcommand{\N}{{\Bbb N}}
\newcommand{\C}{{\Bbb C}}
\newcommand{\Z}{{\Bbb Z}}
\newcommand{\R}{{\Bbb R}}
\newcommand{\Rp}{{\R_+}}
\newcommand{\eps}{\varepsilon}
\newcommand{\supp}{\operatorname{supp}}
\newcommand{\la}{\langle}
\newcommand{\ra}{\rangle}
\newcommand{\const}{\operatorname{const}}
\renewcommand{\emptyset}{\varnothing}
\newcommand{\di}{\partial}
\newcommand{\hotimes}{\hat\otimes}

\renewcommand{\author}[1]{\medskip{\Large #1}\par\medskip}

\newcommand{\pii}{\pi_{\nu\otimes\sigma}}
\newcommand{\RR}{{\cal R}}
\newcommand{\RX}{{\RR\times X}}
\newcommand{\ZZ}{\Z_{+,\,0}^\infty}

\begin{center}{\Large \bf Orthogonal decompositions for L\'evy processes
with\\ an application to the gamma, Pacsal, and Meixner
 processes}\end{center}

\author{Eugene Lytvynov}

\noindent{\sl Institut f\"{u}r Angewandte Mathematik,
Universit\"{a}t Bonn, Wegelerstr.~6, D-53115 Bonn,\\ Germany; SFB
611, Univ.~Bonn, Germany;  BiBoS, Univ.\ Bielefeld, Germany\\[2mm]
 {\rm E-mail: lytvynov@wiener.iam.uni-bonn.de}}



\begin{abstract}

\noindent It is well known  that between all processes with
independent increments, essentially only the Brownian motion  and
the Poisson process possess the chaotic representation property
(CRP). Thus, a natural question appears: What is an appropriate
analog of the CRP in the case of a general L\'evy process. At
least three approaches are possible here. The first one, due to
 It\^o, uses the CRP of the Brownian motion and the
Poisson process, as well as the representation of a L\'evy process
through those processes. The second approach, due to Nualart and
Schoutens, consists in representing any square-integrable random
variable as a sum of multiple stochastic integrals constructed
with respect to a family of orthogonalized centered power jumps
processes. The third approach, never applied before to the L\'evy
processes, uses the idea of orthogonalization of polynomials with
respect to a probability measure defined on the dual of a nuclear
space. The main aims of the present paper are to develop the three
approaches in the case of a general ($\R$-valued) L\'evy process
on a Riemannian manifold and (what is more important) to
understand a relationship between these approaches. We apply the
obtained results to the gamma, Pascal, and Meixner processes, in
which case  the analysis related to the orthogonalized polynomials
becomes essentially simpler and richer than in the general case.

\end{abstract}


\section{Introduction} It is well known (\cite{De, Em, 4}) that between all processes with independent increments,
essentially only the Brownian motion and the Poisson process
possess the chaotic representation property  (CRP) (see e.g.\
\cite{mey} for the definition of the CRP). Thus, in the situation
where one has to deal with a general L\'evy process, a natural
question appears: What is an appropriate analog of the CRP in this
case. At least three approaches are possible here.

The first approach was proposed by It\^o in  \cite{Ito}, see also
\cite{De}. Let ${\bf X}=({\bf X}_t)_{t\ge0}$ be a square
integrable L\'evy process. One may always choose a version of this
process which is right-continuous, has limits from the left, and
does not have fixed discontinuities. By the L\'evy--Khintchine
formula, the process $\bf X$ can be decomposed  into a Gaussian
process plus a stochastic integral with respect to a  compensated
Poisson process on $\R_+\times \R$. By taking the tensor product
of the chaos decomposition of the Brownian and Poisson components,
one obtains a unitary isomorphism between the  canonical
$L^2$-space of the process $\bf X$ and a symmetric Fock over
$L^2(\R_+\times\R;\vartheta)$, where the measure $\vartheta$ is
derived from the L\'evy--Khintchine formula. In what follows, we
will call this approach the Fock space decomposition for the
L\'evy process.

The second approach is due to Nualart and Schoutens \cite{NS} (see
also the recent book \cite{S}) and consists in the following.
Denote  $\tilde\nu(ds){:=}s^{2}\nu(ds)$, where $\nu$ is the L\'evy
measure of the process $\bf X$, and suppose that the Laplace
transform of the measure $\tilde\nu$ may be extended to an
analytic function in a neighborhood of zero, which particularly
implies that the measure $\tilde\nu$ has all moments finite.
Instead of considering a single L\'evy process $\bf X$, one
defines through the original process a family of centered power
jump processes $({\bf X}^{(m)}_t)_{t\ge0}$, $m\in\N$, which are
called by the authors Teugels martingales. (We note that these
processes have also been used in \cite{LyReShch} in a study of
Wick calculus of compound Poisson measures.) One then shows that
there exist processes $$ {\bf Y}^{(m)}_t= {\bf
X}_t^{(m)}+a_{m,\,m-1}{\bf X}^{(m-1)}_t+\dots+a_{m,\, 1}{\bf
X}^{(1)}_t,\qquad t\ge0,\ m\in\N, $$ with numbers $a_{m,\,
n}\in\R$ independent of $t$, which are orthogonal for different
$m$'s. A theorem in \cite{NS} states that any square-integrable
random variable admits a representation as a sum of multiple
stochastic integrals with respect to the processes $({\bf
Y}^{(m)}_t)_{t\ge0}$, $m\in\N$. In the case where the CRP holds,
either all the processes $({\bf X}_t^{(m)})_{t\ge0}$, $m\in\N$,
coincide (Poisson case), or ${\bf X}_t^{(m)}=0$, $t\ge0$, $m\ge2$
(Gaussian case), and thus one arrives at the classical chaos
decomposition. As main examples of application of this approach,
the L\'evy processes of Meixner's type---the gamma, Pascal, and
Meixner processes---were considered in \cite{NS,S}.

It should be noticed that the Pascal and Meixner processes were
originally introduced in \cite{BR} and \cite{ST}, respectively. In
\cite{Gri}, the Meixner process was proposed for a model for risky
asserts and an analog of the Black--Sholes formula was
established.  We also refer to the recent paper \cite{Tsil} and
the references therein, where many  properties of the gamma
process are discussed in detail.

Finally, the third approach  is based on the idea  of
orthogonalization of polynomials with respect to a measure defined
on the dual of a nuclear space. In the case of a general
probability measure, such a procedure was first proposed by
Skorohod in \cite[Sect.~11]{Sko}. Suppose that we are given a
L\'evy noise measure $\mu$ on the space ${\cal D}'={\cal
D}'(\R_+)$ of distributions on $\R_+$.
 We
then have a realization of the L\'evy process given by ${\bf
X}_t{:=}\la\cdot,\chi_{[0,t]}\ra$, $t\ge0$, where $\chi_A$ denotes
the indicator of a set $A$. Suppose that  the L\'evy measure of
the process satisfies the same condition as in the approach of
Nualart and Schoutens. Though the space generated by the multiple
stochastic integrals with respect to the process $\bf X$ is, in
general,  a proper subspace of $L^2({\cal D}';\mu)$, the set of
all continuous polynomials, i.e., finite sums of functions of the
form ${\cal D}'\ni\omega\mapsto\la\omega^{\otimes n},f_n\ra$ with
$f_n$ being a symmetric test function on $\R_+^n$ (denoted by
$f_n\in{\cal D}^{\hotimes n}$), is dense in $L^2({\cal D}';\mu)$.
Denoting by ${\cal P}_n^\sim({\cal D })$ the closure of all
continuous polynomials of power $\le n$ in $L^2({\cal D}';\mu)$,
we obtain the orthogonal decomposition
\begin{equation}\label{21212121} L^2({\cal D}';\mu)=\bigoplus _{n=0}^ \infty
{\bf P}_n({\cal D}'),\end{equation} where ${\bf P}_n({\cal D}')$
denotes the orthogonal difference ${\cal P }^\sim_n({\cal
D}')\ominus {\cal P}^\sim_{n-1}({\cal D}')$.  The set of all
projections ${:}\la\cdot^{\otimes n},f_n\ra{:}$ of continuous
monomials $\la \cdot^{\otimes n},f_n\ra$ onto ${\bf P}_n({\cal D
}')$ is dense in ${\bf P}_n({\cal D}')$. Therefore, we can define,
for each $n\in\N$, a Hilbert space ${\frak F}_n$
 as the closure of the set ${\cal D}^{\hotimes n}$ in the
 norm generated by the scalar product \begin{equation}\label{gffdz} \la
 f_n,g_n\ra_{{\frak F}_n}{:=}\frac1{n!}\, \int_{{\cal D}'}{:}\la \omega^{\otimes n},f_n\ra{:}
 \, {:}\la\omega^{\otimes n},g_n\ra{:}\,\mu(d\omega),\qquad f_n,g_n\in{\cal D}^{\hotimes n}.
 \end{equation} Hence, we
 naturally arrive at a unitary operator between $L^2({\cal D}';\mu)$
and the Hilbert space \begin{equation}\label{qztztzt}{\frak
F}{:=}\bigoplus_{n=0}^\infty{\frak F}_n\,n!,\end{equation} where
${\frak F}_0{:=}\R$. In the case where the CRP holds true, ${\frak
F}$ is the usual symmetric Fock space over $L^2(\R_+;dx)$, while
the obtained unitary operator coincides with the operator given
through the multiple stochastic integral decomposition. In the
general case, the main problem is to identify the scalar product
\eqref{gffdz} explicitly.

As well known, for a probability measure $\tilde \mu$ on $\R$ such
that the set of all polynomials on $\R$ is dense in
$L^2(\R;\tilde\mu)$, the procedure of orthogonalization of
polynonials with respect to $\tilde \mu$ is equivalent to finding
the Jacobi matrix $J$ defining  an (unbounded) self-adjoint
operator in $\ell_2$, whose spectral measure is $\tilde\mu$, see
\cite[Ch.~VII, Sect.~1]{b}. More exactly, the Jacobi matrix
$J=(J_{i,j})_{i,j=0}^\infty$ is given by $J_{i,j}=0$ if $|i-j|>1$,
$J_{i,i-1}=J_{i-1, i}=\beta_i$, $i\in\N$, $J_{i,i}=\alpha_i$,
$i\in\Z_+$, where $\alpha_i,\beta_i$ are the coefficients of the
recurrence relation satisfied by the system $\{Q_i,\ i\in\Z_+\}$
of the normalized orthogonal polynomials:
\begin{equation}\label{ftft9098}
sQ_i(s)=\beta_{i+1}Q_{i+1}(s)+\alpha_i Q_i(s)+\beta_iQ_i(s),\qquad
i\in\Z_+,\ Q_{-1}(s){:=}0,\ Q_1(s){:=}1.\end{equation} Vice versa,
given a Jacobi matrix $J$ which determines a self-adjoint operator
in $\ell_2$, there exists a spectral measure $\tilde\mu$ of $J$,
which is a probability measure on $\R$ whose normalized orthogonal
polynomials satisfy the recurrence relation \eqref{ftft9098} with
the coefficients $\alpha_i,\beta_i$ defined by $J$ as above.

In the infinite-dimensional case, the role of a Jacobi matrix
should be played by a Jacobi field $(A(\varphi))_{\varphi\in{\cal
D }}$,  see the works by Berezansky {\it et al.}\
\cite{bere,berre,BeLi} and \cite{Ly} for the notion of a Jacobi
field. More precisely, each operator $A(\varphi)$ in $\frak F$
should correspond in the functional realization to the operator of
multiplication by the monomial $\la\cdot,\varphi\ra$, and  should
have a three-diagonal structure with respect to the orthogonal
decomposition \eqref{qztztzt}, i.e., each $A(\varphi)$ should be a
sum of a creation, neutral, and annihilation operator. However,
for a general probability measure on the dual of a nuclear space,
the problem of existence of a corresponding Jacobi field is still
open. We refer here to the work \cite{Chap} where a sufficient
condition for the existence of the Jacobi field was given in terms
of the moments of the measure.

This problem was solved for the  gamma process in \cite{KL} (see
also \cite{silva}) and for the processes of Meixner's
type---defined even on a general manifold $X$ (instead of
$\R_+$)---in the  recent paper \cite{Ly3} (see also the paper
\cite{bere3} for the case of the Pascal process on $\R_+$). More
precisely, in \cite{KL,Ly3} the Hilbert space $\frak F$ was {\it a
priori\/} defined as an extended Fock space, a Jacobi field in
${\frak F}$ was constructed, and then it was shown that the
spectral measure of the Jacobi field is a corresponding L\'evy
process.

 In \cite{beme1}, it was shown that the extended Fock
space is naturally isomorphic to a direct sum of subspaces of
$L^2$-spaces of a special form. In \cite{beme2,beme3}, it was
shown that the extended Fock space decomposition of the gamma
process can be thought  of as an expansion of any $L^2$-random
variable in multiple integrals constructed by using the family of
the resolutions of the identity of the operators of the
corresponding Jacobi field.

The main aims of the present paper are to develop the three
approaches in the case of a general ($\R$-valued) L\'evy process
without  Brownian part on a manifold, and (what is more important)
to understand a relationship between these approaches. So, the
contents of the present paper is as follows.

In Section~\ref{uisdgzu}, we present a definition of a L\'evy
process  on a general Riemannian manifold $X$. In the case where
the L\'evy measure of the process has the first moment finite, we
essentially follow the definition and construction of the process
in \cite{Tsil}, using the corresponding Poisson process. In the
case of the infinite first moment (which, for a L\'evy process on
$\R_+$, yields that the trajectories of the process are of
unbounded variation on any finite  interval of time), we define a
L\'evy process as a  generalized process on a   space ${\cal D}'$
of distributions on $X$ (which is dual of a nuclear space) through
its law---a probability measure $\mu$ on ${\cal D}'$ given by its
Fourier transform (compare with \cite[Ch.~III, Sec.~4]{GV}).

In Section~\ref{hkftd}, using the CRP of the  Poisson process on
$\R\times X$ with intensity $\nu\otimes\sigma$, we construct a
unitary operator between the space $L^2({\cal D}';\mu)$ and the
symmetric Fock space over $L^2(\R\times X;\nu\otimes\sigma)$. Here
$\nu$ is the L\'evy measure of the process and $\sigma$ is its
intensity measure.

In Section~\ref{sec3}, using the unitary operator mentioned above,
we prove the Nualart--Schoutens chaotic decomposition for the
L\'evy process on the manifold. We note that even in the standard
case where $X=\R$ or $\R_+$, our proof differs from the original
one in \cite{NS}. Furthermore, we discuss the unitary description
of $L^2({\cal D}';\mu)$ which appears from the obtained chaotic
decomposition (the original description of \cite{NS} works only in
the case of the one-dimensional underlying space). Using our
approach, we, in particular, easily derive a formula for
multiplication of any multiple stochastic integral by a monomial
of first order.

Next, in Section~\ref{sec4}, we derive from the Nualart--Schoutens
decomposition, the decomposition \eqref{21212121} for  $L^2({\cal
D }';\mu)$. More exactly, we explicitly identify the scalar
product \eqref{gffdz}, and furthermore, we write down a
representation of any function ${:}\la \cdot^{\otimes
n},f_n\ra{:}$ as defined above through multiple stochastic
integrals as in Section~\ref{sec3}. Thus, we establish a
correspondence between the second and third approaches to chaotic
decomposition (Corollaries~\ref{is58u} and \ref{hgf}). Our results
also allow one to identify the Jacobi field corresponding to the
L\'evy process.

Finally, in Section~\ref{434343}, we apply the obtained results to
the processes of Meixner's type---the gamma, Pascal, and Meixner
processes. These are characterized by e.g.\ a special form of
their L\'evy measure $\nu$, more exactly, the measure
$\tilde\nu(ds)=s^{2}\,\nu(ds)$ is a probability measure on $\R$
whose orthogonal polynomials with leading coefficient 1 are
polynomials of Meixner's type and satisfy the recurrence relation
\eqref{ghfd} below. In turn, their one-dimensional distributions
happen to be again orthogonality measures of polynomials of
Meixner's type. We show that, for these processes, analysis
related to the orthogonalized polynomials becomes  sufficiently
simpler and reacher than in the general case. (Notice that this
difference cannot be  felt if one restricts himself only to the
approach of It\^o or that of Nualart and Schoutens!) In
particular, we re-derive the Jacobi field of the process as a
special case of the general formulas. As a by-product of our
considerations, we obtain a more explicit description of the
structure of the extended Fock space $\frak F$ than that obtained
in \cite{beme1}. Finally, we show that, in the case of the gamma
and Pascal processes, the decomposition of each function ${:}\la
\cdot^{\otimes n},f_n\ra{:}$ into a sum of multiple stochastic
integrals of the Nualart--Schoutens type may be interpreted as
decomposition of a random measure obtained by dividing the whole
space into (non-random) disjoint parts (Theorem~\ref{cfstrtr}).

\section{L\'evy processes on manifolds}\label{uisdgzu}

In this section, we construct a L\'evy process on a manifold,
using ideas of  \cite{GV,Tsil}.

 Let $X$ be a complete, connected, oriented $C^\infty$
(non-compact) Riemannian manifold and let ${\cal B}(X)$ be the
Borel $\sigma$-algebra on $X$. Let $\sigma$ be a Radon measure on
$(X,{\cal B}(X))$ that is non-atomic, i.e., $\sigma(\{x\})=0$ for
every $x\in X$ and non-degenerate, i.e., $\sigma(O)>0$ for any
open set $O\subset X$. (We note that the assumption of the
nondegeneracy of $\sigma$ is nonessential and the results below
may be generalized to the case of a degenerate $\sigma$.) Note
that $\sigma(\Lambda)<\infty$ for each $\Lambda\in{\cal
B}_c(X)$---the set of all sets from ${\cal B}(X)$ with compact
closure.

Let  ${\cal R}{:=}\R\setminus\{0\}$. We endow $\RR$ with the
relative topology of $\R$ and let ${\cal B}(\RR)$ denote the Borel
$\sigma$-algebra on $\RR$. Let $\nu$ be a Radon measure on
$(\RR,{\cal B}(\RR))$, whose  support contains an infinite number
of points. Let $\tilde\nu(ds){:=}s^2\,\nu(ds)$. We suppose that
$\tilde \nu$ is a finite measure on $(\RR,{\cal B}(\RR))$, and
furthermore,   there exists $\varepsilon>0$
 such that
\begin{equation}\label{4335ew4} \int_{\RR} \exp\big(\varepsilon
|s|\big)\,\tilde\nu(ds)<\infty.\end{equation} By \eqref{4335ew4},
the Laplace transform of the measure $\tilde\nu$ is well defined
in a neighborhood of zero and may be extended to an analytic
function on $\{z\in\C: |z|<\varepsilon\}$. Therefore, the measure
$\tilde\nu$ has all moments finite, and moreover, the set of all
polynomials is  dense in $L^2(\RR;\tilde\nu)$. Next, we evidently
have
 \begin{equation}\label{hjgfcahzf} \forall n\ge2:\qquad \int_\RR |s|^n
 \,\nu(ds)<\infty.\end{equation}   Notice that \eqref{hjgfcahzf}
implies that $\nu(|s|\ge a)<\infty$ for any $a>0$.

We will first additionally suppose that
\begin{equation}\label{ftzedre}\int_\RR |s|\,\nu(ds)<\infty.\end{equation}

Let $\Gamma_\RX$ denote the configuration space over $\RX$ defined
as follows: $$ \Gamma_\RX{:=}\big\{\gamma\subset \RR\times X :\,
\sharp(\gamma\cap\{|s|\ge\varepsilon\}\times\Lambda)<\infty\text{
for each $\varepsilon>0$ and $\Lambda\in{\cal B}_c(X)$}\big\}.$$
Here, $\sharp(A)$ denotes the cardinality of a set $A$. Each
$\gamma\in\Gamma_\RX$ may be identified with the positive Radon
measure $$\sum_{(s,x)\in\gamma}\delta_{(s,x)}\in{\cal M}_+(\RX),$$
where $\delta_{(s,x)}$ denotes the Dirac measure with mass at
$(s,x)$, $\sum_{(s,x)\in\varnothing}\delta_{(s,x)}{:=}$zero
measure, and\linebreak  ${\cal M}_+(\RX)$ denotes the set of all
positive Radon measures on the Borel $\sigma$-algebra ${\cal
B}(\RX)$. We endow the space $\Gamma_\RX$ with the relative
topology as a subset of the space ${\cal M}_+(\RX)$ with the vague
topology, i.e., the weakest topology on $\Gamma_\RX$ with respect
to which all maps
$$\Gamma_\RX\ni\gamma\mapsto\la\gamma,f\ra{:=}\int_\RX
f(s,x)\,\gamma( ds,dx) =\sum_{(s,x)\in\gamma} f(s,x) $$ are
continuous. Here, $f\in C_0(\RX)$(${:=}$the set of all continuous
functions on $\RX$ with compact support).

Let $\pii$ denote the Poisson measure on $(\Gamma_\RX,{\cal
B}(\Gamma_\RX))$  with intensity $\nu\otimes\sigma$. This measure
can be characterized by its Fourier transform
\begin{equation}\label{jsduz} \int_{\Gamma_\RX}
e^{i\la\gamma,f\ra}\,\pii(d\gamma)=\exp\bigg[\int_\RX
(e^{if(s,x)}-1)\,\nu(ds)\,\sigma(dx)\bigg],\qquad f\in C_0(\RX).
\end{equation}
We refer to e.g.\ \cite{Kingman} for an explicit construction of
the Poisson measure.

Since the measure $\sigma$ is non-atomic, it follows from the
construction of the Poisson measure that, for $\pii$-a.e.\
$\gamma\in\Gamma_\RX$, \begin{equation}\label{hjvcfg} \forall
(s,x),\, (s',x')\in\gamma:\ (s,x)\ne(s',x')\Rightarrow x\ne
x'.\end{equation} We fix an arbitrary $x_0\in X$ and let $
B(x_0,r)$ denote the closed ball in $X$ centered at $x_0$ and of
radius $r$. For any $n,k\in \N$, we then have, by the Mecke
identity (e.g.\ \cite{Kingman}),
\begin{align} \int_{\Gamma_\RX}\int_\RX |s|^k
\chi_{B(x_0,n)}(x)\,\gamma(ds,dx)\, \pii(d\gamma)&= \int_\RX |s|^k
\chi_{B(x_0,n)}(x)\,\nu(ds)\,\sigma(dx)\notag\\&= \int_\Rp
|s|^k\,\nu(ds)\,\sigma(B(x_0,n))<\infty.\label{hjvfzgf}\end{align}
We denote by $\widetilde \Gamma_\RX\in{\cal B}(\Gamma_\RX)$ the
set of all $\gamma\in\Gamma_{\RX}$ for which \eqref{hjvcfg} holds
and\linebreak $\int_\RX |s|^k
\chi_{B(x_0,n)}(x)\,\gamma(ds,dx)<\infty$ for $k,n\in \N$. By
\eqref{hjvfzgf}, we get $\pii(\widetilde\Gamma_\RX)=1$. Let ${\cal
B }(\widetilde\Gamma_\RX)$ denote the trace $\sigma$-algebra of
${\cal B}(\Gamma_\RX)$ on $\widetilde\Gamma_\RX$.

For each $\gamma\in \widetilde\Gamma_\RX$, we define $
\omega(\gamma){:=}\sum_{(s,x)\in\gamma}s\,\delta_x$,  which is a
signed  Radon measure on $X$. Furthermore, the mapping
$\widetilde\Gamma_\RX\ni\gamma\mapsto \omega(\gamma)\in{\cal M}(X)
$ is Borel-measurable. Here, ${\cal M}(X)$ denotes the space of
all signed Radon measures on $X$ endowed with the vague topology.
Let $\Omega(X)$ denote the image of $\widetilde\Gamma_\RX$ under
the mapping $\gamma\mapsto\omega(\gamma)$ and let ${\cal
B}(\Omega(X))$ denote the trace $\sigma$-algebra of ${\cal
B}({\cal M}(X))$ on $\Omega(X)$.

We define a L\'evy process on $X$ with intensity measure $\sigma$
and L\'evy measure $\nu$ as a generalized process on $\Omega(X)$
whose law is the probability measure $\mu_{\nu,\,\sigma}$ on
$(\Omega(X),{\cal B}(\Omega(X)))$ obtained as the image of
$\pi_{\nu\otimes\sigma}$ under the measurable mapping
\begin{equation}\label{hjfscgf}
\widetilde\Gamma_{\RX}\ni\gamma\mapsto
\omega(\gamma)\in\Omega(X).\end{equation} As follows from
\eqref{jsduz}, the Fourier transform of $\mu_{\nu,\,\sigma}$ is
given by $$ \int_{\Omega(X)} e^{i\la s,\varphi\ra}\,
\mu_{\nu,\,\sigma}(d\omega)=\exp\bigg[\int_\RX(e^{is\varphi(x)}-1)\,\nu(ds)\,\sigma(dx)\bigg],\qquad
\varphi\in C_0(X).$$ Here, $C_0(X)$ denotes the set of all
continuous functions on $X$ with compact support.

In the case where \eqref{ftzedre} does not hold, such a direct
construction of a L\'evy process is, of course, impossible, so we
proceed as follows.

We denote by ${\cal D}$ the space $C_0^\infty(X)$ of all
real-valued infinite differentiable functions on $X$ with compact
support.  This space may be naturally endowed with a topology of a
nuclear space, see e.g.\ \cite{BUS} for the case $X=\R^d$ and
e.g.\  \cite{Die} for the case of a general Riemannian manifold.
We recall that
\begin{equation}\label{uiguz}{\cal D}=\operatornamewithlimits{proj\,lim}_{\tau\in
T}{\cal H}_\tau.\end{equation} Here, $T$ denotes the set of all
pairs $(\tau_1,\tau_2)$ with $\tau_1\in\Z_+$ and $\tau_2\in
C^\infty(X)$, $\tau_2(x)\ge1$ for all $x\in X$, and ${\cal
H}_{\tau}={\cal H}_{(\tau_1,\tau_2)}$ is the Sobolev space on $X$
of order $\tau_1$ weighted by the function $\tau_2$, i.e., the
scalar product in ${\cal H}_\tau$, denoted by $(\cdot,\cdot)_\tau$
is given by
\begin{equation}\label{hjdrt}(f,g)_{\tau}=\int_X \bigg(
f(x)g(x)+\sum_{i=1}^{\tau_1} \langle
\nabla^{i}f(x),\nabla^{i}g(x)\ra_{T_x(X)^{\otimes{i}}}\bigg)\tau_2(x)\,dx,\end{equation}
where $\nabla^i$ denotes the $i$-th (covariant) gradient, and $dx$
is the volume measure on $X$. For $\tau,\tau'\in T$, we will write
$\tau'\ge\tau$ if $\tau'_1\ge\tau_1$ and  $\tau'_2(x)\ge\tau_2(x)$
for all $x\in X$.

 The space $\cal D$ is densely and continuously
embedded into the real $L^2$-space $L^2(X;\sigma)$. As easily
seen, there always exists $\tau_0\in T$ such that ${\cal
H}_{\tau_0}$ is continuously embedded into $L^2(X;\sigma)$. We
denote $T'{:=}\{\tau\in T: \tau\ge \tau_0\}$ and  \eqref{uiguz}
holds with $T$ replaced by $T'$. In what follows, we will just
write $T$ instead of $T'$. Let ${\cal H}_{-\tau}$ denote the dual
space of ${\cal H}_\tau$ with respect to the zero space
$L^2(X;\sigma)$. Then ${\cal
D}'=\operatornamewithlimits{ind\,lim}_{\tau\in T}{\cal H}_{-\tau}$
is the dual of ${\cal D}$ with respect to $L^2(X;\sigma)$, and we
thus get the standard triple $$ {\cal D}'\supset L^2(X;\sigma)
\supset {\cal D}.$$ The dual pairing between any $\omega\in{\cal
D}'$ and $\xi\in{\cal D}$ will be denoted by $\la\omega,\xi\ra$.
We can evidently consider ${\cal M}(X)$ a subset of ${\cal D}'$ by
identifying any $\omega\in{\cal M}(X)$ with $\tilde\omega\in{\cal
D }'$ by setting $\la\tilde\omega,\varphi\ra{:=}\int_X
\varphi(x)\,\omega(dx)$ for each $\varphi\in{\cal D}$. Let ${\cal
C}({\cal D}')$ denote the cylinder $\sigma$-algebra on $ {\cal
D}'$. Then, the trace $\sigma$-algebra of ${\cal C}({\cal D}')$ on
${\cal M}(X)$ coincides with ${\cal B}({\cal M}(X))$. Thus, any
probability measure on ${\cal B}({\cal M}(X))$ may also be
considered as a probability measure on ${\cal C}({\cal D}')$.

 We now
define a centered L\'evy process as a generalized process  on
${\cal D}'$ whose law is the probability measure on $({\cal
D}',{\cal C}({\cal D}'))$ given by its Fourier transform
 \begin{equation}\label{rew4w}
\int_{{\cal D}'} e^{i\la \omega,\varphi\ra}\,
\rho_{\nu,\,\sigma}(d\omega)=\exp\bigg[\int_\RX(e^{is\varphi(x)}-1-is\varphi(x))\,\nu(ds)\,\sigma(dx)\bigg],\qquad
\varphi\in {\cal D}.\end{equation} The existence of
$\rho_{\nu,\,\sigma}$ follows from the Bochner--Minlos theorem. In
the case where \eqref{ftzedre} holds, $\rho_{\nu,\,\sigma}$
coincides with the measure obtained by centering
$\mu_{\nu,\,\sigma}$.

\section{The Fock space decomposition for a L\'evy process}\label{hkftd}

In this section, we will discuss the Fock space decomposition  for
a L\'evy processes which comes from the multiple stochastic
integral decomposition for the corresponding  Poisson process.

Let first \eqref{ftzedre} hold. Since the mapping \eqref{hjfscgf}
is one-to-one and since $\mu_{\nu,\,\sigma}$ is the image measure
of $\pi_{\nu\otimes\sigma}$ under \eqref{hjfscgf}, we conclude
that the operator $$ L^2(\widetilde\Gamma_\RX;\pi_{\nu\otimes
\sigma})\ni F\mapsto UF=(UF)(\omega){:=}F(\gamma(\omega))\in
L^2(\Omega(X);\mu_{\nu,\,\sigma})$$ is unitary. Here,
$\omega\mapsto \gamma(\omega)$ is the inverse mapping of
\eqref{hjfscgf}.  As well known (see e.g.\ \cite{Sur}), the
Poisson measure $\pi_{\nu\otimes\sigma}$ possesses the chaotic
decomposition property. More exactly, for each $g_n\in
L^2(\RX;\nu\otimes\sigma)^  {\hotimes n}$, $n\in\N$, one can
construct a multiple stochastic integral $I^{(n)}(g_n)$ with
respect to the centered Poisson process on $\RX$ with intensity
$\nu\otimes\sigma$. Here, $\hotimes$ denotes the symmetric tensor
product. Furthermore, we have \begin{equation}\label{tre}
\|I^{(n)}(g_n)\|^2_{L^2(\widetilde\Gamma_\RX;\,\pi_{\nu\otimes\sigma})}=n!\,
\|g_n\|^2_{L^2(\RX;\,\nu\otimes\sigma)^{\hotimes
n}},\end{equation} and any random variable $G\in
L^2(\widetilde\Gamma_\RX;\pi_{\nu\otimes\sigma})$  can be uniquely
represented as a sum of multiple stochastic integrals: $$
G=\sum_{n=0}^\infty I^{(n)}(g_n),$$ where
$I^{(0)}(g_0){:=}g_0\in\R$, and the series converges in
$L^2(\widetilde\Gamma_\RX;\pi_{\nu\otimes\sigma})$. Thus, we have
the unitary operator \begin{equation}\label{zrftzrf} {\cal
F}(L^2(\RX;\nu\otimes\sigma))\ni g=(g_n)_{n=0}^\infty \mapsto
Ig{:=}\sum_{n=0}^\infty I^{(n)}(g_n)\in
L^2(\Gamma_\RX;\pi_{\nu\otimes\sigma}).\end{equation} Here, for a
real Hilbert space ${\cal H}$, ${\cal F}({\cal H})$ denotes the
real Fock space over $\cal H$: $${\cal F}({\cal
H}){:=}\bigoplus_{n=0}^\infty {\cal H}^{\hotimes n}\,n!,\qquad
{\cal H }^{\hotimes 0}{:=}\R$$ (for a Hilbert space $H$ and a
constant $c>0$, we denote by $Hc$ the Hilbert space consisting of
the same elements as $H$ and with scalar product equal to $c$
times the scalar product in $H$).

Setting $J{:=}UI$, we get the unitary operator
\begin{equation}\label{hfutf} {\cal F}(L^2(\RX;\nu\otimes\sigma))\ni
g=(g_n)_{n=0}^\infty \mapsto Jg=\sum_{n=0}^\infty J^{(n)}(g_n)\in
L^2(\Omega(X);\mu_{\nu,\,\sigma}),\end{equation} where $
J^{(n)}(g_n){:=}U(I^{(n)}(g_n))$.

Let us consider the general case. We  still have the unitary
operator $I$ as in \eqref{zrftzrf}. We now denote by ${\cal
P}_{\mathrm cyl}({\cal D}')$ the set of all cylindrical
polynomials on ${\cal D}'$, i.e., the set of all functions on
${\cal D}'$ that are finite sums of constants and monomials of the
form $$ \la\cdot,\varphi_1\ra\dotsm \la\cdot,\varphi_n\ra,\qquad
\varphi_1,\dots,\varphi_n\in{\cal D},\ n\in\N.$$

\begin{lemma}\label{111}
${\cal P}_{\mathrm cyl}({\cal D}')$ is a dense subset of
$L^2({\cal D}';\rho_{\nu,\,\sigma})$\rom.\end{lemma}

\noindent{\it Proof}. Let ${\cal D}_\C$ denote the
complexification of the real space $\cal D$, and define the
Laplace transform of the measure $\rho_{\nu,\,\sigma}$ by
\begin{equation}\label{uzfztfdztdzdffd} L(\theta){:=}\int_{{\cal
D}'}e^{\la\omega,\theta\ra}\,\rho_{\nu,\,\sigma}(d\omega),\qquad
\theta\in{\cal D}_\C,\end{equation} provided the integral on the
right hand side of \eqref{uzfztfdztdzdffd} exists. Using
\eqref{4335ew4} and \eqref{rew4w}, we conclude that the Laplace
transform $L$ is well defined and analytic in a neighborhood of
zero  in ${\cal D}_\C$ (which equivalently means that $L$ is
bounded on this neighborhood and G-holomorphic, see e.g.\
\cite{Dineen}). Then, by (the proof of) \cite[Sec.~10,
Th.~1]{Sko}, we conclude the statement. \quad $\blacksquare$

In what follows, when writing
$I^{(n)}(g_n(s_1,\dots,s_n)f_n(x_1,\dots,x_n))$, we will
understand under $s_1$,\dots, $s_n$, $x_1$,\dots, $x_n$ the
variables in which the integration in the multiple stochastic
integral is carried out.

\begin{lemma}\label{wqwawaw} For any $\varphi_1,\dots,\varphi_n\in{\cal
D}$\rom, $n\in\N$\rom, the distribution of the $\R^n$-valued
random variable $(\la\cdot,\varphi_1\ra,\dots,
\la\cdot,\varphi_n\ra)$ under $\rho_{\nu,\,\sigma}$ coincides with
the distribution of the $\R^n$-valued random variable\linebreak
$(I^{(1)}(s\varphi_1(x)),\dots,I^{(1)}(s\varphi_n(x)))$ under
$\pi_{\nu\otimes\sigma}$\rom.
 \end{lemma}

\noindent {\it Proof}.  By using \eqref{jsduz} and \eqref{rew4w},
we see that, for any $\varphi\in{\cal D}$, the Fourier transform
of the random variable $\la\cdot,\varphi\ra$ under
$\rho_{\nu,\,\sigma}$ coincides with the Fourier transform of the
random variable $I^{(1)}(s\varphi(x))$ under
$\pi_{\nu\otimes\sigma}$.  Therefore, by linearity, we conclude
that, for any fixed $\varphi_1,\dots,\varphi_n\in{\cal D}$,
$n\in\N$, the Fourier transform of the random variable
$(\la\cdot,\varphi_1\ra,\dots, \la\cdot,\varphi_n\ra)$ under
$\rho_{\nu,\,\sigma}$ coincides with the Fourier transform of the
 random variable
$(I^{(1)}(s\varphi_1(x)),\dots,I^{(1)}(s\varphi_n(x)))$ under
$\pi_{\nu\otimes\sigma}$. From here the statement follows.\quad
$\blacksquare$

For any $f\in L^2(\RX;\nu\otimes\sigma)$,  let $A(f)$ denote the
operator in ${\cal F}(L^2(\RX;\nu\otimes\sigma))$ whose image
under the unitary $I$ is the operator of multiplication by the
random variable $I^{(1)}(f)$. We then have (see e.g.\ \cite{Sur})
\begin{equation}\label{dr6se4zwa43z}
A(f)g_n=A^+(f)g_n+A^0(f)g_n+A^-(g)g_n,\qquad n\in\Z_+
\end{equation} where \begin{align}
A^+(f)g_n(s_1,x_1,\dots,s_{n+1},x_{n+1})&=
(f(s_1,x_1)g_n(s_2,x_2,\dots,s_{n+1},x_{n+1}))^\sim,\label{ffdtf}\\
A^0(f)g_n(s_1,x_1,\dots,s_{n},x_{n})&=n
(f(s_1,x_1)g_n(s_1,x_1,\dots,s_{n},x_{n}))^\sim,\label{8zr76}\\
A^-(f)g_n(s_1,x_1,\dots,s_{n-1},x_{n-1})&= n\int_{\RX}
f(s,x)g_n(s,x,s_1,x_1,\dots,s_{n-1},x_{n-1})\,\nu(ds)\,\sigma(dx),\label{ersd4}
\end{align}
that is, $A^+(f)$, $A^0(f)$, $A^-(f)$ are creation, neutral, and
annihilation operators in the Fock space, respectively. Here,
$(\cdot)^\sim$ denotes the symmetrization of a function,
and\linebreak  $g_n\in L^2(\RX;\nu\otimes\sigma)^{\hotimes n}$ is
such  that the right hand side of \eqref{8zr76} belongs
to\linebreak $L^2(\RX;\nu\otimes\sigma)^{\hotimes n}$.

\begin{lemma}\label{drsees} The linear span of the set \begin{equation}\label{cftdxts}  \big\{1,\
I^{(1)}(s\varphi_1(x))\dotsm I^{(1)}(s\varphi_n(x)),\
 \varphi_1,\dots,\varphi_n\in{\cal D},\
n\in\N\big\}\end{equation} is dense  in
$L^2(\widetilde\Gamma_{\RR\times X};\pi_{\nu\otimes\sigma})$\rom.
\end{lemma}

\noindent{\it Proof}. By  \eqref{dr6se4zwa43z}--\eqref{ersd4},
each function $I^{(1)}(s\varphi_1(x))\dotsm
I^{(1)}(s\varphi_n(x))$, $\varphi_1,\dots,\varphi_n\in{\cal D}$,
$n \in\N$, indeed belongs to $L^2(\widetilde\Gamma_{\RR\times
X};\pi_{\nu\otimes\sigma})$.
 Let $\frak L$ denote the closed linear span
of the set \eqref{cftdxts}, and thus we have to show that ${\frak
L }=L^2(\widetilde\Gamma_{\RR\times X};\pi_{\nu\otimes\sigma})$.

Let us consider the unitary operator
 \begin{equation}\label{hgtz} L^2({\cal R};\nu)\ni
f=f(s)\mapsto \frac{f(s)}s\in L^2({\cal
R};\tilde\nu).\end{equation} As we already mentioned above, the
set of functions $\{s^n,\ n\in\Z_+\}$ is total in $L^2({\cal
R};\tilde\nu)$, i.e., its linear span is a dense subset.
Therefore, the set of functions  $\{s^n,\ n\in\N\}$ is total in
$L^2({\cal R};\nu)$. Hence, it  suffices to show that, for each
$f_n\in L^2(X;\sigma)^{\otimes n}$, $n\in\N$, and
$k_1,\dots,k_n\in\N$,
\begin{equation}\label{zfztd} I^{(n)}\big((s_1^{k_1}\dotsm
s_n^{k_n}f_n(x_1,\dots,x_n))^\sim\big)\in {\frak L}.\end{equation}

Let ${\cal O}_c(X)$ denote the algebra of sets in $X$ generated by
all open sets in $X$ with compact closure. For each $m\in\N$, we
introduce a random measure ${\cal X}^{(m)}$ on $X$ by setting, for
each $\Delta\in{\cal O}_c(X)$, $${\cal
X}^{(m)}(\Delta){:=}I^{(1)}(s^m\chi_\Delta(x)).$$ It is easy to
see (see also  the proof of Lemma~\ref{lem2} below) that, for any
$f_n\in L^2(X;\sigma)^{\otimes n}$, $n\in\N$, and any
$k_1,\dots,k_n\in\N$, $$ I^{(n)}\big((s_1^{k_1}\dotsm
s_n^{k_n}f^{(n)}(x_1,\dots,x_n))^\sim\big)=\int_{X^n}
f^{(n)}(x_1,\dots,x_n)\, d{\cal X}^{(k_1)}(x_1)\dotsm d{\cal
X}^{(k_n)}(x_n),$$ where the expression on the right hand side
denotes the multiple stochastic integral constructed with respect
to the random measures ${\cal X}^{(k_1)},\dots,{\cal X}^{(k_n)}$.
Therefore, it is enough to prove that, for any disjoint
$\Delta_1,\dots,\Delta_n\in{\cal O}_c(X)$, $$ {\cal
X}^{(k_1)}(\Delta_1)\dotsm{\cal X}^{(k_n)}(\Delta_n)\in {\frak
L}.$$ This, in turn, will follow from the next

{\it Claim}.  For $m\in\N$\rom, let ${\frak L}_m$ denote the
closure in $L^{2m}(\widetilde\Gamma_{\RR\times
X};\pi_{\nu\otimes\sigma})$ of the linear span of the set
\eqref{cftdxts}. (In particular, we get ${\frak L}={\frak L}_1$.)
Then\rom, we have for any $n,m\in\N$ and $\Delta\in{\cal O}_c(X)$
\begin{equation}\label{gutf} {\cal X}^{(n)}(\Delta)\in{\frak L}_m.\end{equation}

{\it Proof of the Claim}. By  \eqref{dr6se4zwa43z}--\eqref{ersd4},
we see that each element of the set \eqref{cftdxts} indeed belongs
to $L^{2m}(\widetilde\Gamma_{\RR\times
X};\pi_{\nu\otimes\sigma})$. We now prove \eqref{gutf} by
induction in $n\in\N$. Let $\Delta\in{\cal O}_c(X)$ and let $n=1$.
Approximate $\chi_\Delta$ by a sequence $\{\varphi_k,\
k\in\N\}\subset{\cal D}$ such that $\bigcup_{k\in\N}\supp
\varphi_k$ is precompact in $X$,
$|\varphi_k(x)|\le\operatorname{const}<\infty$ for all $k\in\N$
and $x\in X$, and $\varphi_k(x)\to\chi_\Delta(x)$ as $k\to\infty$
for each $x\in X$. We then get by
\eqref{dr6se4zwa43z}--\eqref{ersd4} \begin{gather}
\int_{\widetilde\Gamma_{\RR\times X}}\big( {\cal
X}^{(1)}(\Delta)-I^{(1)}(s\varphi_k(x))\big)^{2m}\,d\pi_{\nu\otimes\sigma}
= \int _{\widetilde\Gamma_{\RR\times X}} I^{(1)}\big(
s(\chi_\Delta(x)-\varphi_k(x))
\big)^{2m}\,d\pi_{\nu\otimes\sigma}\notag\\  =
\big(A^{2m}(s(\chi_\Delta(x)-\varphi_k(x)))\Omega,\Omega\big)_{{\cal
F}(L^2(\RX;\,\nu\otimes\sigma))}\to 0\quad \text{as
}k\to\infty.\label{usasdxztutdf}
\end{gather}
Here, $\Omega{:=}(1,0,0,\dots)$ denotes the vacuum in the Fock
space.

Suppose the statement holds for ${\cal X}^{(1)},\dots,{\cal
X}^{(n)}$ and let us prove it for ${\cal X}^{(n+1)}$. By
\eqref{dr6se4zwa43z}--\eqref{ersd4} $$ {\cal
X}^{(n+1)}(\Delta)={\cal X}^{(1)}(\Delta){\cal
X}^{(n)}(\Delta)-I^{(2)}\big((s_1s_2^n\chi_\Delta(x_1)\chi_\Delta(x_2))^\sim\big)-\int_\RR
s^{n+1}\,\nu(ds)\,\sigma(\Delta).$$ We evidently have ${\cal
X}^{(1)}(\Delta){\cal X}^{(n)}(\Delta)\in{\frak L}_m$ for all
$m\in \N$. Hence, it remains to show that
$$I^{(2)}\big((s_1s_2^n\chi_\Delta(x_1)\chi_\Delta(x_2))^\sim\big)\in
{\frak L}_m,\qquad m\in\N.$$ Clearly,
$$I^{(2)}\big((s_1s_2^n\chi_\Delta(x_1)\chi_\Delta(x_2))^\sim\big)=
I^{(2)}\big((s_1s_2^n\chi_{\Delta^2\setminus D}(x_1,
x_2))^\sim\big),$$ where $D{:=}\{(x_1,x_2)\in X^2: x_1=x_2\}$. We
approximate the indicator $\chi_{\Delta^2\setminus D}(x_1,x_2)$ by
a sequence of functions $\{f_k,\ k\in\N \}$ such that each $f_k$
is a finite sum of functions
$\chi_{\Delta_1}(x_1)\chi_{\Delta_2}(x_2)$ with
$\Delta_1,\Delta_2\in{\cal O}_c(X)$,
$\Delta_1\cap\Delta_2=\varnothing$, $\bigcup_{k\in\N}\supp f_k$ is
precompact in $X^2$, $|f_k(x_1,x_2)|\le1$ for all $k\in\N$ and
$(x_1,x_2)\in X^2$, and $f_k(x_1,x_2)\to\chi_{\Delta^2\setminus
D}(x_1,x_2)$ as $k\to\infty$ for all $(x_1,x_2)\in X^2$. For any
$\Delta_1,\Delta_2\in{\cal O}_c(X)$,
$\Delta_1\cap\Delta_2=\varnothing$, $$
I^{(2)}\big((s_1s_2^n\chi_{\Delta_1}(x_1)\chi_{\Delta_2}(x_2))^\sim\big)={\cal
X }^{(1)}(\Delta_1){\cal X}^{(n)}(\Delta_2)\in{\frak L}_m,\qquad
m\in\N,$$ which yields
\begin{equation}\label{iugzdrdt}I^{(2)}\big((s_1s_2^n f_k(x_1,x_2))^\sim\big) \in{\frak L}_m,
\qquad k, m\in\N,\end{equation} We clearly have
\begin{equation}\label{zftrdr}I^{(2)}\big((s_1s_2^n f_k(x_1,x_2))^\sim\big) \to I^
{(2)}\big((s_1s_2^n\chi_{\Delta^2\setminus D}(x_1,
x_2))^\sim\big)\quad\text{in } L^2(\widetilde \Gamma_{\RR\times X
};\pi_{\nu\otimes\sigma})\end{equation} as $k\to\infty$. By
\eqref{iugzdrdt} and \eqref{zftrdr}, we will get the inclusion
$$I^{(2)}\big((s_1s_2^n\chi_{\Delta^2\setminus D}(x_1,
x_2))^\sim\big)\in {\frak L}_m,\qquad m\in\N,$$ provided we show
that  $\{I^{(2)}\big((s_1s_2^n f_k(x_1,x_2))^\sim\big),\ k\in\N\}$
is a Cauchy sequence in each $L^{2m}(\widetilde \Gamma_{\RR\times
X };\pi_{\nu\otimes\sigma})$, $m\in\N$. But this can be easily
derived, analogously to \eqref{usasdxztutdf}, by using the formula
which expresses a product of arbitrary multiple stochastic
integrals with respect to the centered Poisson process through a
corresponding sum of multiple stochastic integrlas, see e.g.\
\cite[Theorem~3]{LRS} or \cite{Sur}. \quad $\blacksquare$

\begin{theorem}\label{jhhzufzt} We may define a unitary operator $${\cal U}:
L^2(\widetilde\Gamma_{\RX};\pi_{\nu\otimes\sigma})\to L^2({\cal
D}';\rho_{\nu,\,\sigma})$$  by setting
\begin{gather*} {\cal U}1{:=}1,\quad {\cal
U}\big(I^{(1)}(s\varphi_1(x))\dotsm
I^{(1)}(s\varphi_n(x))\big){:=}\la\cdot,\varphi_1\ra\dotsm
\la\cdot,\varphi_n\ra,\\\varphi_1,\dots,\varphi_n\in{\cal D},\
n\in\N. \end{gather*} Furthermore\rom, by setting ${\cal
J}{:=}{\cal U}I$\rom, we get a unitary  operator $$ {\cal
F}(L^2(\RX;\nu\otimes\sigma))\ni g=(g_n)_{n=0}^\infty\mapsto {\cal
J}g=\sum_{n=0}^\infty {\cal J }^{(n)}(g_n)\in L^2({\cal
D}';\rho_{\nu,\,\sigma}),$$ ${\cal J }^{(n)}$ denoting the
restriction of ${\cal J}$ to $L^2(\RX;\nu\otimes\sigma)^{\hotimes
n}$\rom.

\end{theorem}

\noindent{\it Proof}. The theorem trvially follows from
Lemmas~\ref{111}--\ref{drsees}.\quad $\blacksquare$

\begin{remark}\label{rdtrdtzf} \rom{
 Evidently, in the case  where \eqref{ftzedre}
holds, the unitary ${\cal J}$ coincides with the operator $J$ as
in \eqref{hfutf} up to the unitary transformation connected with
the centering of the measure.}\end{remark}

\section{The Nualart--Schoutens chaotic decomposition for a L\'evy
processes}\label{sec3}

In this section, we will generalize the result of Nualart and
Schoutens \cite{NS} (see also \cite[Section~5.4]{S})  concerning a
chaotic decomposition for a usual L\'evy process on $\R$ to the
case of a L\'evy process on the manifold $X$.

We introduce polynomials \begin{equation}\label{gfrtspioiop}
P_n(s)=s^{n}+a_{n,\,n-1}s^{n-1}+a_{n,n-2}s^{n-2}+\dots+a_{n,1}s,\qquad
a_{n,\,i}\in\R,\ i=1,\dots,n-1,\ n\in\N, \end{equation} in such a
way that $$ \int_\RR P_n(s)P_m(s)\, \nu(ds)=0\qquad\text{if }n\ne
m.$$ Using unitary \eqref{hgtz}, we see that
\begin{equation}P_n(s)=\widetilde P_{n-1}(s)s,\label{uzasxfztdfxsz}\end{equation} where $(\widetilde
P_n(\cdot))_{n=0}^\infty$ is the system of polynomials  with
leading coefficient 1 that are orthogonal with respect to the
measure $\tilde\nu(ds)$ on $\RR$. We also evidently have that
$(P_n(\cdot))_{n=1}^\infty$ is a total system in $L^2(\RR;\nu)$.
Then, by Theorem~\ref{jhhzufzt}, the system of the random
variables consisting of the constants ${\cal J}^{(0)}(f_0)$,
$f_0\in \R$, and $$ {\cal J}^{(n)}\big((P_{k_1}(s_1)\dotsm
P_{k_n}(s_n)f_n(x_1,\dots,x_n))^\sim\big),\qquad
k_1,\dots,k_n\in\N,\ f_n\in L^2(X;\sigma)^{\otimes n},\ n\in\N, $$
is total in $L^2({\cal D}';\rho_{\nu,\, \sigma})$.

Denote by $\ZZ $ the set of all sequences $\alpha$ of the form
$\alpha=(\alpha_1,\alpha_2,\dots,\alpha_n,0,0,\dots)$,
$\alpha_i\in\Z_+$, $n\in\N$. Let $|\alpha|{:=}\sum_{i=1}^\infty
\alpha_i$, evidently $|\alpha|\in\Z_+$. For $\alpha\in\ZZ$, denote
$$ P_{\alpha} (s_1,\dots,s_{|\alpha|}){:=}
\underbrace{P_1(s_1)\dotsm P_1(s_{\alpha_1})}_{\text{$\alpha_1$
times}} \,\underbrace{P_2(s_{\alpha_1+1})\dotsm
P_2(s_{\alpha_1+\alpha_2})}_{\text{$\alpha_2$ times}}\dotsm$$ if
$|\alpha|\in\N$, and $P_{\alpha}(s_1,\dots,s_{|\alpha|}){:=}1$ if
$|\alpha|=0$. We  then see that the system of the random variables
\begin{gather*} {\cal I} ^{\alpha}(f_{\alpha}){:=} {\cal
J}^{(|\alpha|)}\big((P_{\alpha}(s_1,\dots,s_{|\alpha|})f_{\alpha}(x_1,\dots,x_{|\alpha|}))^\sim\big)
,\\  f_{\alpha}\in L^2(X;\sigma)^{\otimes |\alpha|},\
\alpha\in\ZZ,\end{gather*} is total in $L^2({\cal
D}';\rho_{\nu,\,\sigma})$. Furthermore, the ${\cal
I}^{\alpha}(\cdot)$'s are pair-wisely orthogonal in $L^2({\cal
D}';\rho_{\nu,\,\sigma})$ for different $\alpha$'s.

Let $$S_n: L^2(\RX;\nu\otimes\sigma) ^{\otimes n}\to L^2_{\mathrm
sym }((\RX)^n;(\nu\otimes\sigma)^{\otimes
n}){:=}L^2(\RX;\nu\otimes\sigma)^{\hotimes n}$$ denote the
symmetrization projection. For $\alpha\in \ZZ$,
$|\alpha|{=:}n\in\N$, denote by
$$L^2_{\alpha}((\RX)^n;(\nu\otimes\sigma)^{\otimes n})$$ the
subspace of  $L^2(\RX;\nu\otimes\sigma) ^{\otimes n}$ consisting
of those functions $g_n(s_1,x_1,\dots,s_n,x_n)$ in
$L^2(\RX;\nu\otimes\sigma) ^{\otimes n}$ which satisfy $$
g_n(s_1,x_1,\dots,s_n,x_n)=g_n(s_{\pi(1)},x_{\pi(1)},\dots,s_{\pi(n)},x_{\pi(n)})$$
for $(\nu\otimes\sigma)^ {\otimes n}$-a.e.\
$(s_1,x_1,\dots,s_n,x_n)\in(\RX)^n$  for any permutation $\pi$ of
$\{1,\dots,n\}$ such that
\begin{equation}\label{hufzdftz}P_{\alpha}(s_1,\dots,s_n)=P_{\alpha}(s_{\pi(1)},\dots,s_{\pi(n)})\qquad
\text{for all }(s_1,\dots,s_n)\in\RR^n.\end{equation} Evidently,
$$L^2_{\alpha}((\RX)^n;(\nu\otimes\sigma)^{\otimes n})=
L^2(\RX;\nu\otimes\sigma)^{\hotimes \alpha_1}\otimes
L^2(\RX;\nu\otimes\sigma)^{\hotimes \alpha_2}\otimes\dotsm\,.$$
Let $$ S_\alpha: L^2(\RX;\nu\otimes\sigma) ^{\otimes n}\to
L^2_\alpha((\RX)^n;(\nu\otimes\sigma)^{\otimes n})$$  denote the
orthogonal projection onto
$L^2_\alpha((\RX)^n;(\nu\otimes\sigma)^{\otimes n})$. Since $$
L^2_{\mathrm sym }((\RX)^n;(\nu\otimes\sigma)^{\otimes n}) \subset
L^2_\alpha((\RX)^n;(\nu\otimes\sigma)^{\otimes n}),$$ we get
$S_n=S_n S_\alpha$. Next, for each $f_n\in L^2(X;\sigma)^{\otimes
n}$, we get by \eqref{hufzdftz} $$ S_\alpha
(P_{\alpha}(s_1,\dots,s_n)f_{n}(x_1,\dots,x_n))=
P_{\alpha}(s_1,\dots,s_n)({\cal S}_\alpha f_{n})(x_1,\dots,x_n),
$$ where ${\cal S}_\alpha$ is the orthogonal projection of
$L^2(X;\sigma)^{\otimes n}$ onto the subspace
$L^2_\alpha(X^n;\sigma^{\otimes n})$ consisting of those functions
$f_{n}\in L^2(X;\sigma)^{\otimes n}$ which satisfy $$
f_{n}(x_1,\dots,x_n)=f_{n}(x_{\pi(1)},\dots,x_{\pi(n)})$$ for
$\sigma^{\otimes n}$-a.e.\ $(x_1,\dots,x_n)\in X^n$ for any
permutation $\pi$ of $\{1,\dots,n\}$ fulfilling \eqref{hufzdftz}.

Consider the operator $$ L^2_\alpha(X^n;\sigma^{\otimes n})\ni
f_{n}\mapsto S_n(P_{\alpha}(s_1,\dots,s_n)f_{n}(x_1,\dots,x_n))\in
L^2_{\mathrm sym}((\RX)^n;(\nu\otimes\sigma)^{\otimes n}). $$ A
direct computation shows that $$
\|S_n(P_{\alpha}(s_1,\dots,s_n)f_{n}(x_1,\dots,x_n))\|^2_{L^2_{\mathrm
sym}((\RX)^n;(\nu\otimes\sigma)^{\otimes
n})}=\frac{\alpha_1!\alpha_2!\dotsm}{n!}\,C_\alpha\,
\|f_{n}\|_{L^2_\alpha(X^n;\,\sigma^{\otimes n})}^2,$$ where
$0!{:=}1$ and \begin{equation}\label{uzedr54t}
C_\alpha{:=}\|P_{1}\|^{2\alpha_1}_{L^2(\RR;\nu)}
\|P_{2}\|^{2\alpha_2}_{L^2(\RR;\nu)}\dotsm.\end{equation}

Thus, by virtue of \eqref{tre} and the definition of ${\cal
I}^\alpha(\cdot)$, we have the following

\begin{lemma}\label{lem1} For each $\alpha\in\ZZ$\rom, $|\alpha|{=:}n\in\N$\rom,  and for each $f_n\in
L^2(X;\sigma)^{\otimes n}$\rom, we have $$ {\cal I}^\alpha
(f_n)={\cal I}^ \alpha({\cal S}_\alpha f_n). $$ Furthermore\rom,
the mapping \begin{equation} L_\alpha^2(X^n;\sigma^{\otimes n})\ni
f_{n}\mapsto {\cal I}^\alpha (f_{n}) \in L^2({\cal
D}';\rho_{\nu,\,\sigma})\label{khacfzt}\end{equation}
 is up to a constant factor an
isometry\rom; more exactly\rom, for each $f_{n}\in
L^2_\alpha(X^n;\sigma^{\otimes n })$ $$ \|{\cal
I}^\alpha(f_{n})\|^2_{L^2({\cal
D}';\,\rho_{\nu,\,\sigma})}=\alpha_1!\,\alpha_2!\dotsm\,
C_\alpha\, \|f_{n}\|_{L_\alpha^2(X^n;\,\sigma^{\otimes n})}^2, $$
where $C_\alpha$ is given by \eqref{uzedr54t}\rom.
\end{lemma}

For each $m\in\N$, we define a random measure ${\bf X}^{(m)}$ on
$X$ by setting, for each $\Delta\in{\cal O}_c(X)$,
\begin{equation}\label{izdr} {\bf X}^{(m)}(\Delta){:=}{\cal U}({\cal
X}^{(m)}(\Delta))={\cal J}^{(1)}(s^m\chi_\Delta
(x)).\end{equation} Notice that, if \eqref{ftzedre} holds, we have
for each $\omega=\sum_{n}s_n \delta_{x_n}\in\Omega(X)$ $${\bf
X}^{(m)}(\Delta,\omega)=\bigg\la \sum_n
s_n^m\delta_{x_n},\chi_\Delta\bigg\ra -\int_\RR
s^m\,\nu(ds)\,\sigma(\Delta).$$ In \cite{NS,S}, in the case
$X=\R$, the  ${\bf X}^{(m)}$, $m\in\N$, was called a Teugels
martingale of order $m$.  Let $${\bf Y}^{(m)}(\Delta){:=}{\cal
J}^{(1)}(P_{m}(s)\chi_\Delta(x)),\qquad \Delta\in{\cal O}_c(X).$$
By \eqref{gfrtspioiop}, $${\bf Y}^{(m)}(\Delta)= {\bf
X}^{(m)}(\Delta)+a_{m,\,m-1}{\bf X}^{(m-1)}(\Delta)+\dots+
a_{m,\,1}{\bf X}^{(1)}(\Delta).$$

\begin{lemma}\label{lem2}
Let $\alpha\in\ZZ$\rom, $|\alpha|{=:}n\in\N$\rom. For each
$f_{n}\in L^2(X^n;\sigma^{\otimes n})$\rom, ${\cal
I}^\alpha(f_{n})$ coincides with the multiple stochastic integral
$$ \int_{X^n} f_{n} (x_1,\dots,x_n)\, d{\bf Y}^{(1)}(x_1)\dotsm
d{\bf Y}^{(1)}(x_{\alpha_1})\,d{\bf Y}^{(2)}(x_{\alpha_1+1})\dotsm
d{\bf Y}^{(2)}(x_{\alpha_1+\alpha_2})\dotsm.$$

\end{lemma}

\noindent  {\it Proof.} For any disjoint
$\Delta_1,\dots,\Delta_n\in{\cal O}_c(X)$, we have by the
definition of a multiple stochastic integral \begin{gather}
 \int_{X^n}\chi_{\Delta_1}(x_1)\dotsm \chi_{\Delta_n}(x_n)
 \, d{\bf
Y}^{(1)}(x_1)\dotsm d{\bf Y}^{(1)}(x_{\alpha_1}) \,d{\bf
Y}^{(2)}(x_{\alpha_1+1})\dotsm d{\bf
Y}^{(2)}(x_{\alpha_1+\alpha_2})\dotsm \notag\\ = {\cal
J}^{(1)}(\Delta_1)\dotsm {\cal J}^{(1)}(\Delta_{\alpha_1}){\cal
J}^{(2)}(\Delta_{\alpha_1+1})\dotsm{\cal
J}^{(2)}(\Delta_{\alpha_1+\alpha_2})\dotsm\,.
\label{one}\end{gather} It follows from Theorem~\ref{jhhzufzt} and
the construction of a multiple stochastic integral with respect to
the centered Poisson process that
\begin{gather*} {\cal U}^{-1} \big({\cal J}^{(1)}(\Delta_1)\dotsm
{\cal J}^{(1)}(\Delta_{\alpha_1}){\cal
J}^{(2)}(\Delta_{\alpha_1+1})\dotsm{\cal
J}^{(2)}(\Delta_{\alpha_1+\alpha_2})\dotsm \big)\\ = I^{(1)}
\big(P_{1}(s)\chi_{\Delta_1}(x)\big)\dotsm
I^{(1)}\big(P_{1}(s)\chi_{\Delta_{\alpha_1}}(x)\big) \\ \times
I^{(1)}\big(P_{2}(s)\chi_{\Delta_{\alpha_1+1}}(x)\big)\dotsm
I^{(1)}\big(P_{2}(s)\chi_{\Delta_{\alpha_1+\alpha_2}}(x)\big)\dotsm\\
= I^{(n)}\big((
P_\alpha(s_1,\dots,s_n)\chi_{\Delta_1}(x_1)\dotsm\chi_{\Delta_n}(x_n)
)^\sim \big).
\end{gather*} Therefore,
\begin{equation}{\cal
J}^{(1)}(\Delta_1)\dotsm {\cal J}^{(1)}(\Delta_{\alpha_1}){\cal
J}^{(2)}(\Delta_{\alpha_1+1})\dotsm{\cal
J}^{(2)}(\Delta_{\alpha_1+\alpha_2})\dotsm = {\cal
I}^{\alpha}(\chi_{\Delta_1}\otimes\dotsm\otimes\chi_{\Delta_n}).\label{two}\end{equation}
By \eqref{one} and \eqref{two}, we have proved the statement for
$f_{n}=\chi_{\Delta_1}\otimes\dots\otimes \chi_{\Delta_n}$. Since
the set of linear combinations of the functions
$\chi_{\Delta_1}\otimes\dots\otimes \chi_{\Delta_n}$ with
$\Delta_1,\dots,\Delta_n\in{\cal O}_c(X)$ disjoint is dense in
$L^2(X;\sigma)^{\otimes n}$, we get the statement by the linearity
and continuity of the mapping $$ L^2(X;\sigma)^{\otimes n}\ni f_n
\mapsto {\cal I}^\alpha (f_n)\in L^2({\cal
D}';\rho_{\nu,\,\sigma})$$ (see Lemma~\ref{lem1}).  \quad
$\blacksquare$

In what follows, we set
$L^2_\alpha(X^{|\alpha|};\sigma^{\otimes|\alpha|}){:=}\R$ for
$|\alpha|=0$.

 \begin{theorem}\label{zgaSZXDT}
Each $F\in L^2({\cal D}';\rho_{\nu,\,\sigma})$ may be uniquely
represented as  a sum of multiple stochastic integrals $$
F=\sum_{\alpha\in\ZZ}{\cal I}^\alpha (f_\alpha),\qquad f_\alpha\in
L^2_\alpha(X^{|\alpha|};\sigma^{\otimes|\alpha|}) ,$$ and $$ {\cal
I}^\alpha(f_\alpha)=\int_{X^{|\alpha|}}f_\alpha(x_1,\dots,x_{|\alpha|})
\, d{\bf Y}^{(1)}(x_1)\dotsm d{\bf Y}^{(1)}(x_{\alpha_1})\, d{\bf
Y}^{(2)}(x_{\alpha_1+1})\dotsm {\bf
Y}^{(2)}(x_{\alpha_1+\alpha_2})\dotsm $$ for $|\alpha|\in\N$ and
${\cal I}^\alpha(f_\alpha)=f_\alpha$ for $|\alpha|=0$\rom.
Furthermore\rom, $$\|F\|^2_{L^2({\cal
D}';\,\rho_{\nu,\,\sigma})}=\sum_{\alpha\in\ZZ}\alpha_1!\,\alpha_2!\dotsm\,
C_\alpha\, \|f_\alpha\|_{L^2(X^n;\,\sigma^{\otimes n})}^2, $$
where $C_\alpha$ is given by \eqref{uzedr54t}\rom.

\end{theorem}

\noindent {\it Proof}. The statement follows by
Theorem~\ref{jhhzufzt} and Lemmas~\ref{lem1}, \ref{lem2}.\quad
$\blacksquare$

We define a Hilbert space \begin{equation}\label{guzzdr} {\bf
H}{:=}\bigoplus_{\alpha\in\ZZ} {\bf H}_\alpha,\qquad  {\bf
H}_\alpha{:=} L^2_\alpha (X^{|\alpha|};\sigma^{\otimes|\alpha|
})\, \alpha_1!\,\alpha_2!\dotsm \, C_\alpha.\end{equation}

As a trivial consequence of  Theorem~\ref{zgaSZXDT}, we get

\begin{corollary}\label{46elkgh} We have the unitary operator $$ {\bf H}\ni
f=(f_\alpha)_{\alpha\in\ZZ}\to {\cal
I}f{:=}\sum_{\alpha\in\ZZ}{\cal I}^\alpha (f_\alpha)\in L^2({\cal
D}';\rho_{\nu,\,\sigma}).$$ \end{corollary}

For each $\varphi\in {\cal D}$, we have $\la
\cdot,\varphi\ra={\cal I}^{(1,0,0,\dots)}(\varphi)$
$\rho_{\nu,\,\sigma}$-a.e., and hence, we define
$\la\cdot,\varphi\ra\in L^2({\cal D}';\rho_{\nu,\,\sigma})$ for
any $\varphi\in L^2(X;\sigma)$ as ${\cal
I}^{(1,0,0,\dots)}(\varphi)$. We will now obtain a formula for the
multiplication of a multiple stochastic integral ${\cal
I}^\alpha(f_\alpha)$, $f_\alpha\in
L^2_\alpha(X^{|\alpha|};\sigma^{\otimes |\alpha|})$, by a random
variable $\la\cdot,\varphi\ra$ (compare with formulas
\eqref{dr6se4zwa43z}--\eqref{ersd4} in the Poisson case).

By the Favard theorem, the system of orthogonal polynomials
$(\widetilde P_n)_{n=0}^\infty$ (see \eqref{uzasxfztdfxsz})
fulfills the recurrence formula \begin{equation*} s\widetilde
P_n(s)=\widetilde P_{n+1}(s)+a_n \widetilde P_n(s)+b_n \widetilde
P_{n-1}(s),\qquad n\in\Z_+,\ \widetilde
P_{-1}(s){:=}0,\end{equation*} with real numbers $a_n$ and
positive numbers $b_n$. Using unitary \eqref{hgtz}, we then get
\begin{equation}\label{zudts} s P_n(s)= P_{n+1}(s)+a_{n-1}
P_n(s)+b_{n-1} P_{n-1}(s),\qquad n\in\N,\
P_{0}(s){:=}0.\end{equation}

For $\alpha\in \ZZ$ and $n\in\N$, we denote $$ \alpha\pm
1_n{:=}(\alpha_1,\dots,\alpha_{n-1},\alpha_n\pm1,\alpha_{n+1},\dots),$$
and let ${\cal I}_\alpha (f_\alpha){:=}0$ if  $\alpha_n<0$ for
some $n\in\N$.

\begin{corollary}\label{tftrtr} Let $\varphi\in L^1(X;\sigma)\cap
L^\infty(X;\sigma)$\rom. Then\rom, for any $\alpha\in\ZZ$\rom,
$|\alpha|{=:}n$\rom, and any $f_\alpha\in
L^2_\alpha(X;\sigma)^{\otimes n }$\rom, we have
\begin{gather}\la\cdot,\varphi\ra \,{\cal I}^\alpha(f_\alpha)={\cal
I}^{\alpha+1_1}\big({\cal
S}_{\alpha+1_1}(\varphi(x_1)f_\alpha(x_2,\dots,x_{n+1}))\big)\notag\\
\text{}+ \frac{\alpha_1!\,(n-\alpha_1)!}{(n-1)!}\, \int_\RR
s^2\,\nu(ds)\, {\cal I}^{\alpha-1_1}\bigg({\cal
S}_{\alpha-1_1}\bigg(\int_X
\varphi(x)f_\alpha(x,x_1,\dots,x_{n-1})\,\sigma(dx)\bigg)\bigg)\notag\\
\text{}+\sum_{n\ge1}\alpha_n\Big[ {\cal
I}^{\alpha-1_n+1_{n+1}}\big( {\cal
S}_{\alpha-1_n+1_{n+1}}(\varphi(x_{\alpha_1+\dots+\alpha_n})f_\alpha(x_1,\dots,x_n)
)\big)\notag\\ \text{} +a_{n-1}{\cal I}^\alpha \big({\cal
S}_\alpha
(\varphi(x_{\alpha_1+\dots+\alpha_n})f_\alpha(x_1,\dots,x_n))\big)
 \notag\\ \text{}+b_{n-1}{\cal I}^{\alpha+1_{n-1}-1_{n}}\big( {\cal
S}_{\alpha+1_{n-1}-1_n}(\varphi(x_{\alpha_1+\dots+\alpha_n})f_\alpha(x_1,\dots,x_n)
)\big)\Big]. \label{uasfzrsf}
\end{gather}

\end{corollary}

\noindent{\it Proof}. The corollary easily follows from the
definition of ${\cal I}^\alpha(\cdot)$, Lemma~\ref{lem1},
\eqref{dr6se4zwa43z}--\eqref{ersd4}, and  \eqref{zudts}.\quad
$\blacksquare$

\section{Orthogonalization of continuous polynomials}\label{sec4}

We denote by ${\cal P}({\cal D}')$ the set of continuous
polynomials on ${\cal D}'$, i.e., functions on ${\cal D}'$ of the
form $$F(\omega)=\sum_{i=0}^n\la\omega^{\otimes
i},f_{i}\ra,\qquad\omega^{\hotimes 0}{:=}1,\ f_{i}\in{\cal
D}^{\hotimes i},\ i=0,\dots,n,\ n\in\Z_+. $$ The greatest number
$i$ for which $f^{(i)}\ne0$ is called the power of a polynomial.
We evidently have ${\cal P}_{\mathrm cyl}({\cal D}')\subset {\cal
P}({\cal D}')$.  We
denote by ${\cal P}_n({\cal D}')$ the set of continuous
polynomials of power $\le n$.

By \cite[Sect.~11]{Sko},  ${\cal P}({\cal D}')$ is a dense subset
of $ L^2({\cal D}';\rho_{\nu\otimes\sigma})$. Let ${\cal P}^\sim
_n({\cal D}')$ denote the closure of ${\cal P}_n({\cal D}')$ in
$L^2({\cal D}';\rho_{\nu,\,\sigma})$,  let ${\bf P}_n({\cal D}')$,
$n\in\N$, denote the orthogonal difference ${\cal P }^\sim_n({\cal
D}')\ominus{\cal P}^\sim_{n-1}({\cal D}')$, and let ${\bf
P}_0({\cal D}'){:=}{\cal P }^\sim_0({\cal D}')$.

\begin{theorem}\label{tew54e} We have the orthogonal
decomposition \begin{equation}\label{tswes} L^2({\cal
D}';\rho_{\nu,\,\sigma})=\bigoplus_{n=0}^\infty{\bf P}_n({\cal
D}'),
\end{equation} and furthermore\rom, \begin{equation}{\bf P}_n({\cal
D}')={\cal I}{\bf H}_n,\label{zhdstrsd}\end{equation} where
\begin{equation}\label{tfzrzzztf} {\bf
H}_n{:=}\bigoplus_{\alpha\in\ZZ:\,
1\alpha_1+2\alpha_2+\cdots=n}{\bf H}_\alpha,\qquad
n\in\Z_+.\end{equation}
\end{theorem}

\noindent{\it Proof}. The orthogonal decomposition \eqref{tswes}
is clear, so we have to prove \eqref{zhdstrsd}, \eqref{tfzrzzztf},
or equivalently
\begin{equation}\label{zudrzd} {\cal P}^\sim_n({\cal
D}')=\bigoplus _{i=0}^n {\cal I}{\bf H}_i =\bigoplus
_{\alpha\in\ZZ: 1\alpha_1+2\alpha_2+\dots\le n}{\cal I}{\bf
H}_\alpha .\end{equation}

 Fix any $\varphi\in{\cal D}$ and consider the continuous
 polynomial $\la \cdot,\varphi\ra^n$. Since ${\cal
 J}^{-1}(\la\cdot,\varphi\ra^n)=A(s\varphi(x))^n\Omega$, we conclude,
 using formulas \eqref{dr6se4zwa43z}--\eqref{ersd4}, that
 \begin{equation}\label{0990zrd}\la\cdot,\varphi\ra^n\in
\bigoplus _{i=0}^n {\cal I}{\bf H}_i.\end{equation} Since the
Laplace transform of the measure $\rho_{\nu,\,\sigma}$ is analytic
in a neighborhood of zero in ${\cal D}_\C$, by
\cite[Lemma~3.9]{KSWY}, there exist $\tau\in T$ and
$\varepsilon_\tau>0$ such that $$ \int_{{\cal D}'}
\exp\big(\varepsilon_\tau \|\omega\|_\tau\big)\,
\rho_{\nu,\,\sigma}(d\omega)<\infty,$$ where $\|\cdot\|_\tau$
denotes the norm in the Hilbert space ${\cal H}_\tau$. Therefore,
 $$ \int_{{\cal D}'}\|\omega\|_\tau^n\,
\rho_{\nu,\,\sigma}(d\omega)<\infty,\qquad n\in\N,$$ which yields
the continuity of the mapping \begin{equation}\label{dsdssd75}
{\cal H}_\tau^{\hotimes n}\ni f_n\mapsto \la \cdot^{\otimes
n},f_n\ra\in L^2({\cal D}';\rho_{\nu,\,\sigma})\end{equation} for
each $n\in\N$. Therefore, any polynomial $\la\cdot^{\otimes
n},f_n\ra$, $f_n\in{\cal D }^{\hotimes n}$, can be approximated in
the $L^2({\cal D}';\rho_{\nu,\,\sigma})$ norm by linear
combinations of polynomials of the form $\la\cdot,\varphi\ra^n$,
$\varphi\in{\cal D }$. Since $\bigoplus _{i=0}^n {\cal I}{\bf
H}_i$ is a linear closed subspace of $L^2({\cal
D}';\rho_{\nu,\,\sigma})$, \eqref{0990zrd} implies the inclusion
${\cal P}^\sim_n({\cal D}')\subset \bigoplus _{i=0}^n {\cal I}{\bf
H}_i$.

Let us prove the inverse inclusion $\bigoplus _{i=0}^n {\cal
I}{\bf H}_i\subset {\cal P}^\sim_n({\cal D}')$. It suffices to
show that, for each $f_m\in L^2(X;\sigma)^{\otimes m}$,
$m\in\{1,\dots,n\}$, $$ {\cal J}^{(m)}\big((s_1^{k_1}\dotsm
s_m^{k_m} f_m(x_1,\dots,x_m))^\sim \big)\in{\cal P}_n^\sim({\cal
D}'),\qquad k_1,\dots,k_m\in\N,\ k_1+\dots+k_m=n.$$ Analogously to
the proof of Lemma~\ref{lem2}, we get $${\cal
J}^{(m)}\big((s_1^{k_1}\dotsm s_m^{k_m} f_m(x_1,\dots,x_m))^\sim
\big) =\int_{X^m}f_m(x_1,\dots,x_m)\,d{\bf X}^{(k_1)}(x_1)\dotsm
d{\bf X}^{(k_m)}(x_m), $$ where the random measures ${\bf
X}^{(i)}$, $i\in\N$, are defined by \eqref{izdr}. Therefore, it is
enough to prove that, for any disjoint
$\Delta_1,\dots,\Delta_m\in{\cal O}_c(X)$,
\begin{equation}\label{ttzl} {\bf X}^{(k_1)}(\Delta_1)\dotsm{\bf
X}^{(k_m)}(\Delta_m)\in{\cal P}^\sim_n({\cal D}'),\qquad
k_1+\dots+k_m=n.\end{equation} But it follows from
Theorem~\ref{jhhzufzt} and  the Claim from the proof of
Lemma~\ref{drsees} that, for any $m,n\in\N$ and any
$\Delta\in{\cal O }_c(X)$,
\begin{equation}\label{jfzd} {\bf X}^{(n)}(\Delta)\in{\cal
P}_n^{\sim\,m}({\cal D}'),\end{equation} where
 ${\cal P}^{\sim\,m} _n({\cal D}')$ denotes the closure
of ${\cal P}_n({\cal D}')$ in $L^{2m}({\cal
D}';\rho_{\nu,\,\sigma})$. Finally, \eqref{jfzd} implies
\eqref{ttzl}.\quad $\blacksquare$

For a monomial $\la\cdot^{\otimes n},f_n\ra$, $f_n\in{\cal
D}^{\hotimes n}$, we denote by ${:}\la \cdot^{\otimes
n},f_n\ra{:}$ the orthogonal projection of $\la\cdot^{\otimes
n},f_n\ra$ onto ${\bf P}_n({\cal D}')$. Since for each monomial
$\la\cdot^{\otimes k},f_k\ra$, $f_k\in{\cal D}^{\hotimes k}$, with
$k<n$, the projection of $\la\cdot^{\otimes k},f_k\ra$ onto ${\bf
P}_n({\cal D}')$ equals zero, the set $\{{:}\la \cdot^{\otimes
n},f_n\ra{:},\ f_n\in{\cal D}^{\hotimes n}\}$ is dense in ${\bf
P}_n({\cal D}')$. Our next aim is to find an explicit
representation of ${:}\la \cdot^{\otimes n},f_n\ra{:}$ through
multiple stochastic integrals ${\cal I}^\alpha(\cdot)$'s.

For each $\alpha\in\ZZ$, $1\alpha_1+2\alpha_2+\dots=n$, $n\in\N$,
and for any function $f_n:X^n\to\R$ we define a function $D_\alpha
f_n:X^{|\alpha|}\to\R$ by setting \begin{align}(D_\alpha
f_n)(x_1,\dots,x_{|\alpha|}){:=}& f(x_1,\dots,x_{\alpha_1},
\underbrace{x_{\alpha_1+1},x_{\alpha_1+1}}_{\text{2 times }},
\underbrace{x_{\alpha_1+2},x_{\alpha_1+2}}_{\text{2 times
}},\dots,
\underbrace{x_{\alpha_1+\alpha_2},x_{\alpha_1+\alpha_2}}_{\text{2
times }},\notag\\ &\quad
\underbrace{x_{\alpha_1+\alpha_2+1},x_{\alpha_1+\alpha_2+1},,x_{\alpha_1+\alpha_2+1}}_{\text{3
times }},\dots).\label{123}\end{align}

\begin{corollary}\label{is58u} For each $f_n\in{\cal D}^{\hotimes
n}$\rom, $n \in\N$\rom, we have \begin{equation}\label{uicvrs}
{:}\la \cdot^{\otimes n},f_n\ra{:}=\sum_{\alpha\in\ZZ:\,
1\alpha_1+2\alpha_2+\dots=n}R_\alpha\, {\cal I}^\alpha (D_\alpha
f_n),\end{equation} where \begin{equation} R_\alpha=
\frac{n!}{\alpha_1!\, (1!)^{\alpha_1}\alpha_2!\,
(2!)^{\alpha_2}\alpha_3!\,
(3!)^{\alpha_3}\dotsm}\,.\label{zufdtrds} \end{equation}
\end{corollary}

\begin{remark}\label{12345}\rom{The $R_\alpha$ given by \eqref{zufdtrds} describes the number of all possible partitions
of a set consisting of $n$ elements into $\alpha_1$ sets
containing one element, $\alpha_2$ sets containing 2 elements, and
so forth. }\end{remark}

\noindent {\it Proof}. Suppose that $f_n=\varphi^{\otimes n}$,
$\varphi\in {\cal D}$. We have \begin{equation}\label{435r8u}
{\cal
U}^{-1}(\la\cdot,\varphi\ra^n)=I^{(1)}(s\varphi(x))^n.\end{equation}
By using \cite[Theorem~2]{LRS} or \cite{Sur}, we express the right
hand side of \eqref{435r8u} as a sum of multiple stochastic
integrals with respect to the centered  Poisson process:
\begin{multline}\label{juhrdtrd} I^{(1)}(s\varphi(x))^n\\ =\sum_{\alpha\in\ZZ:\, 1\alpha_1+2\alpha_2+\dots=n}
R_\alpha I^{(|\alpha|)} \big((s_1 \dotsm s_{\alpha_1}
s^2_{\alpha_1+1}\dotsm s^2_{\alpha_1+\alpha_2}\dotsm (D_\alpha
\varphi^{\otimes n}
)(x_1,\dots,x_{|\alpha|}))^\sim\big)+G_\alpha,\end{multline} where
\begin{equation}\label{rdrdrd}G_\alpha\in \bigoplus _{\alpha\in\ZZ :\,
1\alpha_1+2\alpha_2+\dots\le n-1} {\cal U}^{-1}{\cal I}{\bf
H}_\alpha.\end{equation} By \eqref{gfrtspioiop}, $P_k(s)-s^k$ is a
polynomial of order $k-1$, which yields by \eqref{juhrdtrd},
\begin{equation} \label{da3232}I^{(1)}(s\varphi(x))^n
=\sum_{\alpha\in\ZZ:\, 1\alpha_1+2\alpha_2+\dots=n} R_\alpha
I^{(|\alpha|)} \big(( P_\alpha(s_1,\dots,s_{|\alpha|})(D_\alpha
\varphi^{\otimes n}
)(x_1,\dots,x_{|\alpha|}))^\sim\big)+G'_\alpha,\end{equation}
where $G_\alpha'$ also belongs to the space from  \eqref{rdrdrd}.
By \eqref{435r8u}, \eqref{da3232}, and Theorem~\ref{tew54e}, we
get \eqref{uicvrs} for $f_n=\varphi^{\otimes n}$.

Next, by the continuity of the mapping \eqref{dsdssd75}, we
conclude that each  mapping \begin{equation}\label{redstar}{\cal
H}_\tau^{\hotimes n}\ni f_n\mapsto {:}\la \cdot^{\otimes
n},f_n\ra{:}\in L^2({\cal D}';\rho_{\nu,\,\sigma})\end{equation}
is also continuous. Let us fix $n\in\N$. Without loss of
generality, we can suppose that, for all $\alpha\in\ZZ$,
$1\alpha_1+2\alpha_2+\dots=n$, the mapping $${\cal H
}_\tau^{\hotimes n}\ni f_n\mapsto D_\alpha f_n\in
L^2(X;\sigma)^{\otimes |\alpha|}$$ is continuous. Now, the formula
\eqref{uicvrs} in the general case  follows by an approximation of
$f_n\in{\cal D}^{\hotimes n}$ by linear combinations of functions
of the form $\varphi^{\otimes n}$, $\varphi\in{\cal D}$, in the
${\cal H}_\tau^{\hotimes n}$  norm.  \quad $\blacksquare$

\begin{corollary}\label{gzuzfd} For any $f_n,g_n\in{\cal D}^{\hotimes n}$\rom,
$n\in\N$\rom, we have \begin{gather}\int_{{\cal
D}'}{:}\la\omega^{\otimes n},f_n\ra{:}\,{:}\la\omega^{\otimes
n},g_n\ra{:}\, \rho_{\nu,\,\sigma}(d\omega)\notag\\
=n!\sum_{\alpha\in\ZZ:\, 1\alpha_1+2\alpha_2+\dots=n}K_\alpha
\int_{X^{|\alpha|}}(D_\alpha f_n)(x_1,\dots,x_{|\alpha|})\notag
\\ \times (D_\alpha g_n)(x_1,\dots,x_{|\alpha|}) \,\sigma^{\otimes
|\alpha|}(dx_1,\dots,dx_{|\alpha|}),\label{tfztftfz}\end{gather}
where \begin{equation}\label{t7rr5r5} K_\alpha=
\frac{n!}{\alpha_1!\,\alpha_2!\dotsm}\,\prod_{k\ge1}\bigg(\frac{\|P_k\|_{L^2(\RR;\nu)}}{k!}\bigg)^{2\alpha_k}.
\end{equation}

\end{corollary}

\noindent{\it Proof}. The corollary follows directly from
Theorem~\ref{zgaSZXDT} and Corollary \ref{is58u}.\quad
$\blacksquare$

Analogously to \eqref{tfzrzzztf}, we define, for each $n\in\Z_+$,
a Hilbert space \begin{equation}\label{jhjhhj}{\frak F
}_n{:=}\bigoplus_{\alpha\in\ZZ:\,
1\alpha_1+2\alpha_2+\dots=n}{\frak F}_{n,\,\alpha},\qquad {\frak
F}_{n,\,\alpha}{:=}
L^2_\alpha(X^{|\alpha|};\sigma^{\otimes|\alpha|})\,K_\alpha.\end{equation}
For each $n\in\N$, we define a mapping $${\cal D}^{\hotimes n}\ni
f_n\mapsto {\cal E}_n f_n\in{\frak F}_n$$ by setting $$ ({\cal
E}_n f_n)_\alpha{:=} D_\alpha f_n,\qquad \alpha\in\ZZ:\,
1\alpha_1+2\alpha_2+\dots=n,$$ where $(\cdot)_\alpha$ denotes the
${\frak F}_{n,\,\alpha}$-component of an element of ${\frak F}_n$.
Since $\{{:}\la\cdot^{\hotimes n}, f_n\ra{:},\ f_n\in{\cal
D}^{\hotimes n }\}$, is a dense subset of ${\bf P}_n({\cal D}')$,
we get by Theorem~\ref{tew54e} and Corollaries~\ref{is58u} and
\ref{gzuzfd} that ${\cal E}_n{\cal D}^{\hotimes n}$ is a dense
subset of ${\frak F}_n$. In what follows, for simplicity of
notations, we will just write $f_n$ instead of $ {\cal E}_nf_n$,
considering, in this way, ${\cal D}^{\hotimes n}$ as a dense
subset of ${\frak F}_n$ (evidently $f_n=g_n$ in ${\cal
D}^{\hotimes n}$ if and only if $f_n=g_n$ in ${\frak F}_n$).

Let $${\frak F}{:=}\bigoplus_{n=0}^\infty {\frak F}_n \,n!,\qquad
{\frak F}_0{:=}\R.$$ There evidently exists a unitary operator
${\cal A}:{\bf H}\to{\frak F}$  which acts between any ${\bf
H}_\alpha$- and ${\frak F}_\alpha$-component of the space ${\bf
H}$, respectively ${\frak F}$ as a constant (depending on
$\alpha$) times the identity operator.

We evidently have the following

\begin{corollary}\label{hgf} We have the unitary operator $$ {\frak U}:{\frak
F}\to L^2({\cal D}';\rho_{\nu,\,\sigma})$$ that is defined through
$$ {\frak U}f_n{:=}{:}\la\cdot^{\otimes n},f_n\ra{:},\qquad
f_n\in{\cal D}^{\hotimes n},\, n\in\Z_+,$$ and then extended by
linearity and continuity to the whole space ${\frak F}$\rom.
Furthermore\rom, $${\frak U}={\cal I}{\cal A}^{-1}.$$

\end{corollary}

\begin{remark}\label{ghgfrdrd}\rom{Using  Theorem~\ref{tew54e} and Corollaries~\ref{tftrtr},
\ref{is58u} one can explicitly identify the Jacobi filed of the
measure $\rho_{\nu,\,\sigma}$, that is, the family of commuting
selfadjoint operators $a(\varphi)$, $\varphi\in{\cal D}$, acting
in the Hilbert space $\frak F$ and satisfying $${\frak
U}a(\varphi){\frak U}^{-1}=\la\cdot,\varphi\ra\cdot,\qquad
\varphi\in{\cal D}, $$ where $\la\cdot,\varphi\ra\cdot$ denotes
the operator of multiplication by $\la\cdot,\varphi\ra$ in
$L^2({\cal D}';\rho_{\nu,\,\sigma})$. As easily seen, for each
$f_n\in{\cal D}^{\hotimes n}$, we have $$
a(\varphi)f_n=a^+(\varphi)f_n+a^0(\varphi)f_n+a^-(\varphi)f_n, $$
where $a^+(\varphi)f_n\in{\frak F}_{n+1}$,
$a^0(\varphi)f_n\in{\frak F}_{n}$, and $a^-(\varphi)f_n\in{\frak
F}_{n-1}$. Though we always have $a^+(\varphi)f_n=\varphi\hotimes
f_n$, i.e., $a^+$ is a usual creation operator, the structure of
the neutral operator $a^0(\varphi)$ and the annihilation operator
$a^-(\varphi)$ is, in general, quite complicated. In the next
section, we will consider a family of L\'evy processes for which
these operators have a much simpler form. }\end{remark}

\begin{remark}\rom{ Let us suppose that \eqref{ftzedre} holds.
For each $n\in\N$, we now define the random measure ${\bf
Y}^{(n)}$ on $X$ by setting, for each $\Delta\in{\cal O}_c(X)$, $$
{\bf Y}^{(n)}(\Delta){:=}J^{(1)}(P_n(s)\chi_\Delta(x)).$$ Then, as
easily seen, all the results of Sections~\ref{sec3} and \ref{sec4}
remain true if we change $L^2({\cal D}';{\rho}_{\nu,\,\sigma})$
for the space $L^2(\Omega(X);\mu_{\nu,\,\sigma})$ and, in the
formulation of Corollary~\ref{tftrtr} and Remark~\ref{ghgfrdrd},
$\la\cdot,\varphi\ra$ for ${:}\la\cdot,\varphi\ra{:}$.
}\label{udtrsrea}\end{remark}

\section{Processes of Meixner's type}\label{434343}

In this section, we consider an application of the above results
to the processes of Meixner's type.

For each $\lambda\ge0$, let $\tilde\nu_\lambda$ denote the
probability measure on $(\R,{\cal B}(\R))$ whose orthogonal
polynomials $(\widetilde P_{\lambda,\, n})_{n=0}^\infty$ with
leading coefficient 1 satisfy the recurrence relation
\begin{gather}\label{ghfd}s\widetilde P_{\lambda,\, n}(s)=
\widetilde P_{\lambda,\, n+1}(s)+\lambda(n+1)\widetilde
P_{\lambda,\, n}(s)+n(n+1)\widetilde P_{\lambda,\, n-1}(s),\\
n\in\Z_+,\, \widetilde P_{\lambda,\,-1}(s){:=}0.\notag\end{gather}
By \cite{Meixner} (see also \cite[Ch.~VI, sect.~3]{Chihara}),
$(\widetilde P_{\lambda,\, n})_{n=0}^\infty$ is a system of
polynomials of Meixner's type, the measure $\widetilde\nu_\lambda$
is uniquely determined by the above condition and is given as
follows: For $\lambda\in[0,2)$, $$ \tilde\nu_\lambda(ds)=
\frac{\sqrt{4-\lambda^2}}{2\pi} \times \big|\Gamma\big(1+i
(4-\lambda^2)^{-1/2}s\big)\big|^2\,\exp\big[-s2(4-\lambda^2)^{-1/2}\arctan
\big(\lambda(4-\lambda^2)^{-1/2}\big) \big]\,ds $$
($\tilde\nu_\lambda$ is a Meixner distribution), for $\lambda=2$
$$\tilde\nu_2(ds)=\chi_{(0,\infty)} (s)e^{-s}s\,ds$$
($\tilde\nu_2$ is a gamma distribution), and for $\lambda>2$
$$\tilde\nu_\lambda(ds)=(\lambda^2-4)\sum_{k=1}^\infty
p_\lambda^k\,k\,\delta_{\sqrt{\lambda^2-4}\,k},\qquad
p_\lambda{:=}\frac{\lambda-\sqrt{\lambda^2-4}}{\lambda+\sqrt{\lambda^2-4}}$$
($\tilde\nu_\lambda$ is now a Pascal distribution).

Since $\tilde\nu_\lambda(\{0\})=0$ for each $\lambda\ge0$, we can
define a L\'evy measure $\nu_\lambda$ by setting
$$\nu_\lambda(ds){:=}\frac{1}{s^2}\,\tilde\nu_\lambda(ds),\qquad
\lambda\ge0.$$ The integral $\int_{\RR} |s|\,\nu_\lambda(ds)$ is
finite for $\lambda\ge 2$ and infinite for $\lambda\in[0,2)$.
Therefore, we will consider, for each $\lambda\ge2$, the L\'evy
process with law $\varrho_\lambda{:=}\mu_{\nu_\lambda,\, \sigma}$,
and for $\lambda\in[0,2)$ the centered L\'evy process with law
$\varrho_\lambda{:=}\rho_{\nu_\lambda,\, \sigma}$.

It follows from \eqref{ghfd} that, for each $\lambda\ge0$, $$
\|\widetilde
P_{\lambda,\,n}\|^2_{L^2(\RR;\,\tilde\nu_\lambda)}=n!\,
(n+1)!\,,\qquad n\in\Z_+, $$ and hence
\begin{equation}\label{fzrd8998} \|P_{\lambda,\,n
}\|^2_{L^2(\RR;\,\nu_\lambda)}=(n-1)!\, n!\,,\qquad
n\in\N,\end{equation} where $P_{\lambda,\,n}(s){:=}\widetilde
P_{\lambda,\,n-1}(s)s$. By \eqref{t7rr5r5} and \eqref{fzrd8998},
we conclude that, for each $\alpha\in\ZZ$,
$1\alpha_1+2\alpha_2+\dots=n$, \begin{equation}\label{bugzu}
K_\alpha=\frac{n!}{\alpha_1!\, 1^{\alpha_1}\alpha_2!\,
2^{\alpha_2}\dotsm}\,.\end{equation}

\begin{remark}\rom{Let us give a combinatoric interpretation of the number $K_\alpha$ in \eqref{bugzu}.
Under a loop $\kappa$ connecting points $x_1,\dots,x_m$,
$m\ge2$, we understand a class of ordered sets
$(x_{\pi(1)},\dots,x_{\pi(m)})$, where $\pi$ is a permutation of
$\{1,\dots,m\}$, which coincide up to a cyclic permutation. 
Let us  also interpret a set $\{x\}$ as a ``one-point'' loop
$\kappa$, i.e., a loop that comes out of $x$. 
Let $\vartheta_n=\{\kappa_1,\dots,\kappa_{|\vartheta_n|}\}$ be a
collection of $|\vartheta_n|$ loops $\kappa_j$ that connect points
from the set $\{x_1,\dots,x_n\}$ so that every point
$x_i\in\{x_1,\dots,x_n\}$
goes into one loop $\kappa_j=\kappa_{j(i)}$ from $\vartheta_n$. 
Then, for $\alpha\in\ZZ$, $1\alpha_1+2\alpha_2+\dots=n$,
$K_\alpha$ is the number of all different collections of loops
connecting points from the set $\{x_1,\dots,x_n\}$ and containing
$\alpha_1$ one-point loops, $\alpha_2$ two-point loops, etc.
}\end{remark}

In the following proposition, we will explicitly identify the
Jacobi field $a_\lambda(\varphi)$, $\varphi\in{\cal D}$, of the
measure $\varrho_\lambda$ (see Remarks~\ref{ghgfrdrd},
\ref{udtrsrea}).

\begin{proposition} For each $\lambda\ge0$\rom, we have for all $\varphi\in{\cal
D}$ and all  $f_n\in{\cal D}^{\hotimes n}$\rom, $n\in\Z_+$\rom,
$$a_\lambda(\varphi)f_n=a^+(\varphi)f_n+\lambda
a^0(\varphi)f_n+a^-(\varphi)f_n.$$ Here\rom, $a^+(\xi)$ is the
standard creation operator\rom: $$
a^+(\varphi)f_n{:=}\varphi\hotimes f_n,$$ $a^0(\varphi)$ is the
standard neutral operator\rom: $$
\big(a^0(\varphi)f_n\big)(x_1,\dots,x_n)=\big(\varphi(x_1)+\dots+\varphi(x_n)\big)f_n(x_1,\dots,x_n),\qquad
n\in\N,\ a^{0}(\varphi)f_0=0, $$ and
$a^-(\varphi)=a^-_1(\varphi)+a^-_2(\varphi)$\rom, where
$a^-_1(\varphi)$ is the standard annihilation operator\rom: $$
\big(a^-_1(\varphi)f_n\big)(x_1,\dots,x_{n-1})=n\int_X \varphi(x)
f_n(x,x_1,\dots,x_{n-1})\,\sigma(dx),\qquad n\in\N,\
a^{0}(\varphi)f_0=0,$$ and  \begin{gather*}
\big(a^-_2(\varphi)f_n\big)(x_1,\dots,x_{n-1})=n(n-1)\big(\varphi(x_1)
f_n(x_1,x_1,x_2,x_3,\dots,x_{n-1})\big)^{\sim} ,\qquad n\ge2,\\
a_2^{-}(\varphi)f_i=0,\quad i=0,1.\end{gather*} \label{323232}
\end{proposition}

\noindent {\it Proof}. The proposition easily follows in the way
described in Remark~\ref{ghgfrdrd} (see also
Remark~\ref{udtrsrea}) through formula \eqref{ghfd}.\quad
$\blacksquare$

\begin{corollary}\label{drtdrd} For each $\lambda\ge0$ and $\omega\in{\cal D}'$\rom,
define $\wick n\in{\cal D}^{'\hotimes n}$ by the recurrence
formula
\begin{gather}\wick{(n+1)}=\wick{n+1}(x_1,\dots,x_{n+1})=\big(\wick
n(x_1,\dots,x_n)\omega(x_{n+1})\big)^\sim\notag\\ \text{}
-n\big(\wick{(n-1)}(x_1,\dots,x_{n-1})\delta(x_{n+1}-x_n)\big)^\sim
\notag\\ \text{} -n(n-1)\big(\wick
{(n-1)}(x_1,\dots,x_{n-1})\delta(x_n-x_{n-1})\delta(x_{n+1}-x_n)\big)^\sim\notag\\
\text{} -\lambda n\big(\wick n
(x_1,\dots,x_n)\delta(x_{n+1}-x_n)\big)^\sim-c_\lambda \big(\wick
n(x_1,\dots,x_n)1(x_{n+1})\big)^\sim,  \notag\\ \wick 0 =1,\ \wick
1=\omega-c_\lambda,\label{udsH} \end{gather} where
$c_\lambda{:=}0$ for $\lambda\in[0,2)$ and
$c_\lambda{:=}2/(\lambda+\sqrt{\lambda^2-4})$ for
$\lambda\ge2$\rom. Then\rom, for each $f_n\in {\cal D}^{\hotimes n
}$\rom, $n\in\N$\rom, $$ {:}\la\omega^{\otimes
n},f_n\ra{:}=\la\wick
n,f_n\ra\qquad\text{\rom{$\varrho_\lambda$-a.e.\ $ \omega\in{\cal
D }'$}}.$$

\end{corollary}

\noindent {\it Proof}. Since $\int_\RR s\nu_\lambda(ds)=c_\lambda$
for each $\lambda\ge2$, the statement trivially holds for $n=1$.
Suppose the statement holds for all $n\le m$ and let us prove it
for $n=m+1$. By Proposition~\ref{323232}, we then get for any
sequence $\{\varphi_i,\,i\in\N\}\subset{\cal D}$
\begin{equation}\label{57r65}{:}\la\omega^{\otimes (m+1)},\varphi_i^{\otimes
(m+1)}\ra{:}=\la\wick {(m+1)},\varphi_i^{\otimes (m+1) }\ra\qquad
i\in\N,\ \text{\rom{$\varrho_\lambda$-a.e.\ $ \omega\in{\cal D
}'$}}.\end{equation} There exist $\tau_1,\,\tau_2\in T$,
$\tau_2>\tau_1$, such that $\varrho_\lambda({\cal H}_{\tau_1})=1$,
for each $n\in\{2,\dots,m+1\}$ the mapping \eqref{redstar} with
$\tau=\tau_2$ is continuous, and for each $\omega\in{\cal
H}_{-\tau_1}$ $\wick{(m+1)}\in{\cal H
}_{-\tau_2}^{\hotimes(m+1)}$. Choose
$\{\varphi_i,\,i\in\N\}\subset{\cal D}$ which is a total set in
${\cal H}_{\tau_2}$. Approximate an arbitrary $f_{m+1}\in{\cal
D}^{\hotimes (m+1)}$ in the ${\cal H}_{\tau_2}^{\hotimes(m+1)}$
topology by  a sequence $\{f_{m+1}^{(k)},\, k\in\N\}$ such that
each $f_{m+1}^{(k)}$ is a linear combination of functions
$\varphi_i^{\otimes(m+1)}$.  Then,  \eqref{57r65} implies the
statement.\quad $\blacksquare$

By using \eqref{udsH} and  \cite{Meixner}, the following
proposition was proved in \cite{Ly3} (see also \cite{KL} and
\cite[Chs.~4, 5]{S}).

\begin{proposition}\label{zfztdz}  The Fourier transform of the
measure $\varrho_\lambda$ is given, in a neighborhood of zero, by
the following formula\rom: for  $\lambda=2$
\begin{equation*}\label{jasuetfzcw}\int_{{\cal
D}'}e^{i\la\omega,\varphi\ra}\,d\mu_2(\omega)=\exp\left[-\int_X
\log(1-i\varphi(x))\, \sigma(dx)\right],\qquad \varphi\in{\cal
D},\ \|\varphi\|_\infty{:=}\sup_{x\in
X}|\varphi(x)|<1,\end{equation*} and for $\lambda\ne2$
\begin{equation*}\label{hisav} \int_{{\cal
D}'}e^{i\la\omega,\varphi\ra}\,d\mu_\lambda(\omega)=\exp\left[-\frac1{\alpha\beta}\int_X\log\bigg(
\frac{\alpha e^{-i\beta\varphi(x)}-\beta
e^{-i\alpha\varphi(x)}}{\alpha-\beta}\bigg)\,\sigma(dx)+ic_\lambda\int_X\varphi(x)\,\sigma(dx)\right]
\end{equation*} for all $\varphi\in{\cal D}$ satisfying \begin{equation*} \label{dyuivf}\bigg\| \frac{\alpha
(e^{-i\beta\varphi}-1)-\beta(e^{-i\alpha\varphi}-1)
}{\alpha-\beta}\bigg\|_\infty<1.\end{equation*} Here\rom, $\alpha,
\beta\in\C$ are defined through the equation $1+\lambda
z+z^2=(1-\alpha z)(1-\beta z)$\rom, $z\in\R$\rom. Furthermore\rom,
we have\rom, for $\lambda=2$\rom,
\begin{equation}\label{ioutz}\sum_{n=0}^\infty \frac1{n!}\,\la {:}\omega^{\otimes n}{:}_2,\varphi^{\otimes
n}\ra=\exp\bigg[-\la
\log(1+\varphi)\ra+\bigg\la\omega,\frac\varphi{\varphi+1}\bigg\ra\bigg],\end{equation}
and for $\lambda\ne 2$
\begin{align}
\sum_{n=0}^\infty \frac1{n!}\,\la\wick n,\varphi^{\otimes
n}\ra&=\exp\bigg[-\frac1{\alpha-\beta}
\bigg\la\log\bigg(\frac{(1-\beta\varphi)^{1/\beta}}{(1-\alpha\varphi)^{1/\alpha}}\bigg)\bigg\ra
\notag\\&\qquad\qquad
\text{}+\frac1{\alpha-\beta}\bigg\la\omega-c_\lambda,\log\bigg(\frac{1-\beta\varphi}{1-\alpha\varphi}\bigg)\bigg\ra
\bigg].\label{hjdtvgr}
\end{align}
Formulas \eqref{hjdtvgr}\rom, \eqref{ioutz} hold for each
$\omega\in{\cal M}(X)$ and for each $\varphi\in{\cal D}$
satisfying $\|\varphi\|_\infty<1$ for \eqref{ioutz} and
$\|\varphi\|_\infty<(\max(|\alpha|,|\beta|))^{-1}$ for
\eqref{hjdtvgr}\rom. More generally\rom, for each fixed $\tau\in
T$\rom, there exists a neighborhood  of zero in ${\cal D}$
\rom(depending on $\lambda$\rom)\rom, denoted by ${\cal
O}_\tau$\rom, such that \eqref{hjdtvgr}\rom, respectively
\eqref{ioutz}, holds for all $\omega\in{\cal H}_{-\tau}$ and all
$\varphi\in  {\cal O}_\tau$\rom.
\end{proposition}

As a direct corollary of this proposition, we get, for each
$\Delta\in{\cal O}_c(X)$, an explicit formula for the distribution
of the random variable $\la\cdot,\chi_\Delta\ra$
 under $\varrho_\lambda$.

\begin{corollary}[\cite{Ly3}] For each $\Delta\in{\cal O}_c(X)$, the
distribution $\varrho_{\lambda,\Delta}$ of the random variable
$\la\cdot,\chi_\Delta\ra$ under $\varrho_\lambda$ is given as
follows\rom: For $\lambda>2$\rom, $\varrho_{\lambda,\Delta}$ is
the negative binomial \rom(Pascal\rom) distribution $$
\varrho_{\lambda,\Delta}=(1-p_{\lambda})^{\sigma(\Delta)}\sum_{k=0}^\infty
 \frac{\big(\sigma(\Delta)\big)_k}{k!}\, p_\lambda^k \, \delta_{\sqrt{\lambda^2-4}\,k}, $$ where
  for $r> 0$ $(r)_0{:=}1$, $(r)_k{:=}r(r+1)\dotsm (r+k-1)$\rom,
$k\in\N$\rom. For $\lambda=2$, $\varrho_{2,\Delta}$ is the Gamma
distribution
$$\varrho_{2,\Delta}(ds)=\frac{s^{\sigma(\Delta)-1}e^{-s}}{\Gamma(\sigma(\Delta))}\,\chi_{(0,\infty)}(s)\,
ds.$$ Finally\rom, for $\lambda\in[0,2)$\rom,
\begin{multline*}\varrho_{\lambda,\Delta}(ds)=\frac{(4-\lambda^2)^{(\sigma(\Delta)-1)/2}}{2\pi
\Gamma(\sigma(\Delta)) }\big|
\Gamma\big(\sigma(\Delta)/2+i(4-\lambda^2)^{-1/2}(s+\lambda\sigma(\Delta)/2\big)
\big|^2\\ \times\exp\big[
-(2s+\lambda\sigma(\Delta))(4-\lambda^2)^{-1/2}\arctan\big(\lambda(4-\lambda^2)^{-1/2}\big)
\big]\,ds.\end{multline*}
\end{corollary}

Let us dwell upon a representation of the random measure ${\bf Y
}^{(n)}$, $n\in\N$, through $\wick n$. We fix $\Delta\in{\cal
O}_c(X)$. We take a sequence $\{\varphi_k,\, k\in\N\}\subset{\cal
D}$ such that the $\varphi_k$'s are uniformly bounded,
$\bigcup_{k\in\N}\supp\varphi_k$ is a precompact set in $X$,
$\varphi_k(x)\to\chi_\Delta(x)$ as $k\to\infty$ for each $x\in X$,
and a sequence $\{\psi_k,\, k\in \N\}\subset {\cal D}^{\otimes n}$
such that $\psi_k$'s are uniformly bounded,
$\bigcup_{k\in\N}\supp\psi_k$ is a precompact set in $X^n$, and
the $\psi_k$'s converge point-wisely to the indicator of the set
$\{(x_1,\dots,x_n)\in B(x_0,r)^n:\ x_1=x_2=\dots= x_n\}$ as
$k\to\infty$. Here, $B(x_0,r)$ is a ball in $X$ such that
$\Delta\subset B(x_0,r)$. We set
$G_k^\sim(x_1,\dots,x_n){:=}{\big(\varphi_k(x_1)\psi_k(x_1,\dots,x_n)\big)^\sim}$.
Evidently,
 $$
G_k^\sim(x_1,\dots,x_n)\to\chi_\Delta(x_1)\chi_{\{x_1=x_2=\dots=x_n\}}(x_1,x_2\dots,x_n)\qquad
\text{for each }(x_1,\dots,x_n)\in X^n.$$ Then, by
Corollary~\ref{is58u} and the majorized convergence theorem, we
easily get $$ {\bf Y}^{(n)}(\Delta)=\lim_{k\to\infty} \la
{:}\cdot^{\otimes n}{:}_\lambda,G_k^\sim\ra \quad \text{in
}L^2({\cal D}';\varrho_\lambda),$$ and we informally write$$ {\bf
Y}^{(n)}(\Delta,\omega)=\la \wick n
(x_1,\dots,x_n),\chi_\Delta(x_1)\chi_{\{x_1=x_2=\dots=x_n\}}(x_1,x_2,\dots,x_n)\ra.$$
Furthermore, by \eqref{udsH}, $\wick n\in{\cal M}(X^n)$ for each
$\omega\in \Omega(X)$. Therefore, in the case where
$\varrho_\lambda$ is concentrated on $\Omega(X)$, i.e.,
$\lambda\ge2$, using the majorized convergence theorem, we
conclude from the above the following

\begin{proposition}\label{677667}
For each $\vartheta_n\in{\cal M}(X^n)$\rom, $n\in\N$\rom, define
$D_n\vartheta_n\in{\cal M}(X)$ by setting\rom, for each
$\Delta\in{\cal B}_c(X)$ $$
D_n\vartheta_n(\Delta){:=}\vartheta_n\big(\{ (x_1,\dots,x_m)\in
X^n:\ x_1\in\Delta,\ x_1=x_2=\dots=x_n \}\big).$$ Let $\lambda\ge
2$\rom. Then\rom, for each $\Delta\in{\cal O}_c(X)$ and
$n\in\N$\rom, $$ {\bf Y}^{(n)}(\Delta,\omega)= {:}\omega{:}_{n,\,
\lambda} (\Delta)\qquad\text{for $\rho_\lambda$-a\rom.e\rom.\
$\omega\in\Omega(X)$},$$ where $$ {:}\omega{:}_{n,\, \lambda}{:=}
D_n\wick n .$$

\end{proposition}

\begin{remark}\rom{Let $\lambda\ge2$ and $\omega\in\Omega(X)$. As easily  seen from \eqref{udsH},  the
${:}\omega{:}_{n,\, \lambda}$'s satisfy the recurrence relation
\begin{gather*} {:}\omega{:}_{n+1,\,
\lambda}=D_2({:}\omega{:}_{n,\,\lambda}\otimes\omega)-\lambda n
\,{:} \omega {:}_{n,\,\lambda}- n^2
\,{:}\omega{:}_{n-1,\,\lambda},\qquad n\in\N,\\
{:}\omega{:}_{0,\,\lambda}=\sigma,\
{:}\omega{:}_{1,\,\lambda}=\omega-c_\lambda \sigma,\end{gather*}
which is, of course, equivalent to the recurrence relation
satisfied by the polynomials $(P_{\lambda,\,n})_{n=0}^\infty$.
}\end{remark}

Finally, we will give a representation of the kernel of a multiple
stochastic integral ${\cal I}^\alpha(f_\alpha)$ through $\wick n$,
where , $n=\alpha_1+2\alpha_2+\cdots$. Just as in
Proposition~\ref{677667}, we suppose that $\lambda\ge2$.

So, let us fix any $\alpha\in\ZZ$, $|\alpha|\in\N$, and let
$n{:=}1\alpha_1+2\alpha_2+\cdots$. Let \begin{multline*}
X_{\alpha}{:=} \big\{\, (x_1,\dots,x_n)\in X^n:\
x_{\alpha_1+1}=x_{\alpha_1+2},\dots,
x_{\alpha_1+2\alpha_2-1}=x_{\alpha_1+2\alpha_2},\\
x_{\alpha_1+2\alpha_2+1}=x_{\alpha_1+2\alpha_2+2}=x_{\alpha_1+2\alpha_2+3},
\dots,
x_{\alpha_1+2\alpha_2+3\alpha_3-2}=x_{\alpha_1+2\alpha_2+3\alpha_3-1}=
x_{\alpha_1+2\alpha_2+3\alpha_3},\\ \dots, x_1\ne x_2\ne\dots\ne
x_{\alpha_1}\ne x_{\alpha_1+2}\ne\dots\ne
x_{\alpha_1+2\alpha_2}\ne x_{\alpha_1+2\alpha_2+3}\ne\dots\ne
x_{\alpha_1+2\alpha_2+3\alpha_3}\ne\cdots\,\big\}.
\end{multline*}
Here, the writing $y_1\ne y_2\ne\dots \ne y_m$ means that $y_i\ne
y_j$ if $i\ne j$. For any permutation $\pi$ of $\{1,\dots,n\}$,
denote $$ X_{\alpha}^\pi{:=}\big\{\,(x_{\pi(1)},\dots,x_{\pi(n)}):
(x_{1},\dots, x_{n})\in X_{\alpha} \,\big\}.$$ Evidently, all the
sets $X_{\alpha}^\pi$ either coincide or disjoint, and there are
exactly $R_\alpha$ disjoint sets between them (see formula
\eqref{zufdtrds} and Remark~\ref{12345}). We denote these disjoint
sets by $X_{\alpha}^{(1)},\dots, X_{\alpha}^{(R_\alpha)}$, where
$X_{\alpha}^{(1)}= X_{\alpha}$.

Now, we fix $n\in\N$ and take all $\alpha\in\ZZ$ such that
$1\alpha_1+2\alpha_2+\cdots=n$. Then,  we get the following
representation of $X^n$ as a union of pair-wisely disjoint sets
$X_{\alpha}^{(i)}$:  \begin{equation}\label{545434w5}
X^n=\bigcup_{\alpha\in\ZZ:\,
1\alpha_1+2\alpha_2+\dots=n}\,\bigcup_{i=1}^{R_\alpha}
X_{\alpha}^{(i)}.\end{equation} Thus, each measure $\vartheta\in
{\cal M}(X^n)$ can be represented as $$ \vartheta
=\sum_{\alpha\in\ZZ:\,
1\alpha_1+2\alpha_2+\dots=n}\,\sum_{i=1}^{R_\alpha}
\vartheta_{\alpha}^{(i)},$$ where
$$\vartheta_{\alpha}^{(i)}(\Delta){:=}\vartheta (\Delta\cap
X^{(i)}_{\alpha}),\qquad \Delta\in{\cal B}(X^n).$$ Let us suppose
that the measure $\vartheta$ is invariant under the action of the
group of permutations of $\{1,\dots,n\}$ on $X^n$, i.e., for each
permutation $ \pi$,  $\vartheta(\Delta)=\vartheta(\pi\Delta)$ for
all $\Delta\in {\cal B}(X^n)$, where
$$\pi\Delta{:=}\{(x_{\pi(1)},\dots,x_{\pi(n)}):(x_1,\dots,x_n)\in\Delta\}.$$
 For any
$f_n\in{\cal D}^{\hotimes n}$, because of symmetricity, we then
obtain: \begin{equation}\label{uhuzuf} \la \vartheta,f_n\ra
=\sum_{\alpha\in\ZZ:\, 1\alpha_1+2\alpha_2+\dots=n}R_\alpha \la
\vartheta_{\alpha},f_n\ra,\end{equation} where
$\vartheta_{\alpha}{:=}\vartheta_{\alpha}^{(1)}$.

Denote $$ \widetilde X^{|\alpha|}{:=}\big\{\,
(x_1,\dots,x_{|\alpha|})\in X^{|\alpha|}:\ x_i\ne x_j\text{ if
}i\ne j\,\big\}.$$ There exits a natural bijection between
$\widetilde X^{|\alpha|}$ and $X_{\alpha}$ given by
\begin{multline*} \widetilde X^{|\alpha|} \ni (x_1,\dots,x_n)\mapsto
T_\alpha
(x_1,\dots,x_n){:=}(x_1,\dots,x_{\alpha_1},\underbrace{x_{\alpha_1+1},x_{\alpha_1+1},}_{\text{2
times}}\dots,\\
\underbrace{x_{\alpha_1+\alpha_2},x_{\alpha_1+\alpha_2},}_{\text{2
times}}
\underbrace{x_{\alpha_1+\alpha_2+1},x_{\alpha_1+\alpha_2+1},
x_{\alpha_1+\alpha_2+1},}_{\text{3 times}}\dots )\in
X_{\alpha}.\end{multline*} Evidently $T_\alpha$ and its inverse,
$T_\alpha^{-1}$, are measurable mappings with respect to the
corresponding trace  $\sigma$-algebras.

Let us consider $\vartheta_{\alpha}$ as a measure on $X_{\alpha}$,
and let $\tilde \vartheta_{\alpha}$ denote the image measure of
$\vartheta_{\alpha}$ under $T_\alpha^{-1}$. As easily seen, for
$f_n\in{\cal D}^{\hotimes n}$,
\begin{equation}\label{zftzrdf} \la \vartheta_{\alpha},f_n\ra
=\la \tilde \vartheta _{\alpha},D_\alpha f_n\ra,\end{equation}
where $D_\alpha$ was defined in \eqref{123}. By \eqref{uhuzuf} and
\eqref{zftzrdf}, taking  $\vartheta= \wick n$ for any $\omega\in
\Omega(X)$, we get
\begin{equation}\label{uuzztzuf} \la \wick n,f_n\ra
=\sum_{\alpha\in\ZZ:\, 1\alpha_1+2\alpha_2+\dots=n}R_\alpha \la
{:}\omega^{\otimes n}{:}_{\lambda,\,\alpha},D_\alpha
f_n\ra,\end{equation} where ${:}\omega^{\otimes
n}{:}_{\lambda,\,\alpha}$ is the corresponding $\tilde
\vartheta_{\alpha}$ measure for $\vartheta=\wick n$.

\begin{theorem}\label{cfstrtr} Let $\lambda\ge 2$, let $\alpha\in\ZZ$\rom,
$m{:=}|\alpha|\in\N$\rom, and let
$n{:=}1\alpha_1+2\alpha_2+\cdots$\rom. Let $f_\alpha\in{\cal
D}^{\hotimes\alpha_1}\otimes {\cal
D}^{\hotimes\alpha_2}\otimes\dotsm\subset {\cal D}^{\otimes
m}$\rom. Then\rom, $${\cal I}^\alpha(f_\alpha)=\la
{:}\cdot^{\otimes
n}{:}_{\lambda,\,\alpha},f_\alpha\ra\qquad\text{$\varrho_\lambda$-a.s.,}$$
where for $\omega\in\Omega(X)$ ${:}\omega^{\otimes
n}{:}_{\lambda,\,\alpha}$ was extended to a measure on $X^{m}$ by
setting it zero on $X^{m}\setminus\widetilde X^{m}$\rom.
\end{theorem}

\begin{remark}\label{gzftdsts}\rom{ It follows from Theorem~\ref{cfstrtr} that the
representation of ${:}\la \omega^{\otimes n},f_n\ra{:}$ obtained
in Corollary~\ref{is58u} in the case of $\varrho_\lambda$,
$\lambda\ge2$, may be understood as taking the partition
\eqref{545434w5} of the space $X^n$, constructing the
corresponding decomposition of each measure $\wick n$ for
$\omega\in\Omega(X)$, and then applying this decomposition to
$\la\wick n,f_n\ra={:}\la \omega^{\otimes n},f_n\ra{:}$\,.
}\end{remark}

\noindent {\it Proof}. We fix $f_\alpha$ as in the formulation of
the theorem.  For any $x_0\in X$, we choose $r>0$ such that $\supp
f_\alpha\subset B(x_0,r)^n$. Let
$\{\varphi_i,\,i\in\N\}\subset{\cal D}^{\otimes m}$ be such that
$|\varphi_i(x_1,\dots,x_m)|\le 1$ for all $i\in\N$ and
$(x_1,\dots,x_m)\in X^m$, $\supp\varphi_i\subset B(x_0,r+ 1)^m$,
$i\in\N$, and $\varphi_i$'s converge point-wisely to the indicator
of the set $\{(x_1,\dots,x_m)\in B(x_0,r)^m: x_1\ne x_2\ne\dots\ne
x_m\}$. For $k\in\{2,\dots,n\}$, let
$\{\psi_i^{(k)},\,i\in\N\}\subset{\cal D}^{\otimes k}$ be such
that $|\psi_i^{(k)}(x_1,\dots,x_k)|\le1$ for all $i\in\N$ and
$(x_1,\dots,x_k)\in X^k$, $\supp \psi_i^{(k)}\subset
B(x_0,r+1)^k$, $i\in\N$, and $\psi_i^{(k)}$'s converge
point-wisely to the indicator of the set $\{(x_1,\dots,x_k)\in
B(x_0,r)^k: x_1=x_2=\dots=x_k\}$. For $i\in\N$, we define
$G_i\in{\cal D}^{\otimes n}$ setting \begin{multline*}
G_i(x_1,\dots,x_n){:=} \big(f_\alpha \varphi_i \big)
(x_1,x_2,\dots,x_{\alpha_1},
x_{\alpha_1+2},x_{\alpha_1+4},\dots,x_{\alpha_1+2\alpha_2},\dots)\\
\times \psi^{(2)}(x_{\alpha_1+1},x_{\alpha_1+2})\dotsm
\psi^{(2)}(x_{\alpha_1+2\alpha_2-1},x_{\alpha_1+2\alpha_2})\\
\times
\psi^{(3)}(x_{1+2\alpha_2+1},x_{1+2\alpha_2+2},x_{1+2\alpha_2+3})\dotsm
\psi^{(3)}(x_{1+2\alpha_2+3\alpha_3-2},x_{1+2\alpha_2+3\alpha_3-1},x_{1+2\alpha_2+3\alpha_3})\dotsm.
\end{multline*}
Let $G^\sim_i\in{\cal D}^{\hotimes n}$ denotes the symmetrization
of $G_i$, $i\in\N$. As easily seen, we get for each
$\omega\in\Omega(X)$:
\begin{equation}\label{2121211} \la \wick n,G_i^\sim\ra\to \la
{:}\omega^{\otimes
n}{:}_{\lambda,\,\alpha},f_\alpha\ra\end{equation} as
$i\to\infty$. On the other hand, by Theorem~\ref{tew54e} and
Corollary~\ref{is58u}, \begin{equation}\label{21211} \la
{:}\cdot^{\otimes n}{:}_\lambda,G_i^\sim\ra\to {\cal
I}^\alpha\big((D_\alpha f_\alpha)\chi_{\{x_1\ne x_2\ne\dots\ne
x_m\}}\big)={\cal I}^\alpha(D_\alpha f_\alpha)\end{equation} as
$i\to\infty$ in $L^2(\Omega(X);\rho_\lambda)$. Comparing
\eqref{2121211} and \eqref{21211}, we conclude the stament.\quad
$\blacksquare$


\noindent {\bf Acknowledgements.} I am  grateful to
Yu.~Berezansky, Yu.~Kondratiev, and D. Mierzejewski for many
useful discussions. I am  also indebted to N. Kachanovsky  for his
remarks on a preliminary version of the paper. The financial
support of the SFB 611 and  the DFG Research Project 436 RUS
113/593 is gratefully acknowledged.

\end{document}